# Siegel measures

By William A. Veech*

## 0. Introduction

The goals of this paper are first to describe and then to apply an ergodic-theoretic generalization of the Siegel integral formula from the geometry of numbers. The general formula will be seen to serve both as a guide and as a tool for questions concerning the distribution, in senses to be made precise, of the set of closed leaves of measured foliations subordinate to meromorphic quadratic differentials on closed Riemann surfaces.

In preparation of a discussion of the main results we recall two earlier theorems. The first of these, by H. Masur, has been a starting point for the present work. Let $q$ be a meromorphic quadratic differential with at worst simple poles on a closed Riemann surface $X$. For a certain countable set of $\theta \in \mathbf{R}$ the horizontal foliation associated to $e^{-2i\theta}q$ has one or more maximal cylinders of closed leaves. Each cylinder determines a pair of vectors $v = \pm re^{i\theta}$, where $r$ is the common $|q|$-length of closed leaves in the cylinder. Let $\Pi(q)$ be the set of vectors, with multiplicities, which arise from closed cylinders as $\theta$ varies. Finally, let $N(q, R) = \mathrm{Card}\{v \in \Pi(q) \mid |v| < R\}$ be the growth function of $\Pi(q)$.

Theorem 0.1 (H. Masur [13], [14]). *Let $(X, q)$ and $N(q, R)$ be as above. There exist $0 < c_1 < c_2 < \infty$ such that*

$$(0.2) \qquad c_1 < \frac{N(q, R)}{R^2} < c_2 \qquad (R \gg 0).$$

In certain instances one can say more with regard to (0.2). Let $G = \mathrm{SL}(2, \mathbf{R})$, and let $\mathcal{U}_q$ be the atlas of natural parameters for $q$ on $X \setminus q^{-1}\{0, \infty\}$. If $g \in G$, an atlas $g\mathcal{U}_q$ is defined by postcomposition of $\mathcal{U}_q$ chart functions with the $\mathbf{R}$-linear transformation $g$. This atlas extends to $X$ as a complex structure and determines a new quadratic differential with the same pattern of zeros and poles as $q$. Denote the new pair by $(X(g), q(g))$, and define $\Lambda(q) \subseteq G$ by $\Lambda(q) = \{g \in G \mid (X, q) \cong (X(g), q(g))\}$, where $\cong$ is biholomorphism identifying $q$ and $q(g)$. $\Lambda(q)$ is a discrete subgroup which is not cocompact but which may be a lattice ([23]). Observe that $\Pi(q(g)) = g\Pi(q)$, $g \in G$.

---

*Research supported by NSF and Universite d'Aix Marseille 2.

THEOREM 0.3 ([23]). *With notations as above if $\Lambda(q)$ is a lattice, there exists $c < \infty$ such that for all $g \in G$*

$$(0.4) \qquad \lim_{R \to \infty} \frac{N(q(g), R)}{R^2} = c\pi.$$

The results of the present paper represent a middle ground of sorts between the general Tchebychev theorem of Masur and the restricted prime geodesic theorem of Theorem 0.3. We consider ergodic actions of $G = \mathrm{SL}(2, \mathbf{R})$ on probability spaces $(X, \mathcal{B}, \mu)$ such that the phase space $X$ is a moduli space of quadratic differentials of norm 1. There is a natural map which assigns to $x \in X$ the set $\Pi(x) \subseteq \mathbf{C}$ associated to (the quadratic differential) $x$ as above. The action is such that $\Pi(gx) = g\Pi(x)$, $g \in G$. We shall prove

THEOREM 0.5. *Let $G$ and $(X, \mathcal{B}, \mu)$ be as above. There exists a constant $c(\mu) < \infty$ such that the following three statements are true*:

I. *Let $\psi \geq 0$ be a Borel function on $\mathbf{R}^2$, and define $\hat{\psi}(x) = \sum_{v \in \Pi(x)} \psi(v)$. Then $\hat{\psi}$ is $\mathcal{B}$-measurable and*

$$(0.6) \qquad \int_X \hat{\psi}(x)\mu(dx) = c(\mu) \int_{\mathbf{R}^2} \psi(u) du.$$

II. *Let $N(x, R)$ be the growth function of $\Pi(x)$. Then*

$$(0.7) \qquad \lim_{R \to \infty} \left\| \frac{N(x, R)}{R^2} - c(\mu)\pi \right\|_1 = 0.$$

III. *If $\psi \in C_c(\mathbf{R}^2)$, then*

$$(0.8) \qquad \lim_{R \to \infty} \left\| \frac{1}{R^2} \sum_{v \in \Pi(x)} \psi\left(\frac{v}{R}\right) - c(\mu) \int_{\mathbf{R}^2} \psi(u) du \right\|_1 = 0.$$

Let $(X, \mathcal{B}, \mu)$ be as in the theorem. One consequence of $L^1(\mu)$-convergence in Parts II–III is the existence of a fixed sequence $R_n \to \infty$ such that the relations (0.7)–(0.8) hold pointwise a.e. when the limits are taken along the sequence $\{R_n\}$. In fact the relation $\Pi(gy) = g\Pi(y)$, $g \in G$, and the $G$-invariance of $\mu$ will imply that for $\mu$-a.e. $y$ the relations (0.7)–(0.8) hold for *all $x = gy$* when the limits are taken along $\{R_n\}$ (Theorem 10.8). In this regard we observe that a countable set $\Pi \subseteq \mathbf{C}$ may have asymptotic growth $cR^2$ without its images $g\Pi$, $g \in G$, having such asymptotic growth, much less with the same constant $c$.

Part I of Theorem 0.5 is reminiscent of and motivated by a classical theorem of Siegel ([20]). Let $G_N = \mathrm{SL}(N, \mathbf{R})$ and $\Gamma_N = \mathrm{SL}(N, \mathbf{Z})$. Equip $G_N/\Gamma_N$ with its normalized Haar measure $\mu_N$. If $\psi \geq 0$ is a Borel function on $\mathbf{R}^N$, and if $e_N$ is (say) the $N^{\mathrm{th}}$ standard basis vector, define $\hat{\psi}$ and $\hat{\psi}_p$ ($p$ for 'primitive') by

$$\hat{\psi}(g\Gamma_N) = \sum_{v \in \mathbf{Z}^N \setminus \{0\}} \psi(gv)$$

$$\hat{\psi}_p(g\Gamma_N) = \sum_{v \in \Gamma_N e_N} \psi(gv).$$

According to Siegel

(0.9) $$\int_{G_N/\Gamma_N} \hat{\psi}(g\Gamma_N) \mu_N(dg\Gamma_N) = \int_{\mathbf{R}^N} \psi(u) du$$

(0.10) $$\int_{G_N/\Gamma_N} \hat{\psi}_p(g\Gamma_N) \mu_N(dg\Gamma_N) = \frac{1}{\zeta(N)} \int_{\mathbf{R}^N} \psi(u) du.$$

To place the Siegel theorem in the context of the present work define $\mathcal{M}_N$ to be the set of Borel measures $\nu$ on $\mathbf{R}^N$ such that $M(\nu) < \infty$, where, setting $N_\nu(R) = \nu(B(0, R))$,

(0.11) $$M(\nu) = \sup_{R>0} \frac{N_\nu(R)}{R^N}.$$

$G_N$ acts naturally by homeomorphisms on $\mathcal{M}_N$ when $\mathcal{M}_N$ is endowed with the $C_c(\mathbf{R}^N)$ weak-$*$ topology. A Borel probability measure $\mu$ on $\mathcal{M}_N$ shall be called a *Siegel measure* if $\mu$ is invariant and ergodic for the $G_N$-action. If $\psi \geq 0$ is a Borel function, define $\hat{\psi}$ on $\mathcal{M}_N$ by duality

$$\hat{\psi}(\nu) = \int_{\mathbf{R}^N} \psi(u) \nu(du).$$

The main theorem for Siegel measures is

THEOREM 0.12. *If $\mu$ is a Siegel measure, there exists a constant $c(\mu) < \infty$ such that*

I. *If $\psi \geq 0$ is a Borel function, then*

(0.13) $$\int_{\mathcal{M}_N} \hat{\psi}(\nu) \mu(d\nu) = c(\mu) \int_{\mathbf{R}^N} \psi(u) du.$$

II. *If $\sigma_N$ is the area of the unit sphere in $\mathbf{R}^N$, then*

(0.14) $$\lim_{R \to \infty} \left\| \frac{N_\nu(R)}{R^N} - c(\mu) \frac{\sigma_N}{N} \right\|_1 = 0.$$

III. *If $\mu$ is supported on $\mathcal{M}_N^e = \{\nu \in \mathcal{M}_N \mid \nu(-U) = \nu(U), U \text{ Borel}\}$, then for all $\psi \in C_c(\mathbf{R}^N)$*

$$(0.15) \quad \lim_{R \to \infty} \left\| \frac{1}{R^N} \int_{\mathbf{R}^N} \psi\left(\frac{v}{R}\right) \nu(dv) - c(\mu) \int_{\mathbf{R}^N} \psi(u) du \right\|_1 = 0.$$

*If $N > 2$, and if $\mu$ is such that $\hat{\psi} \in L^2(\mu)$ for all $\psi \in C_c(\mathbf{R}^N)$, then convergence in (0.14)–(0.15) also holds pointwise a.e. $\mu$.*

Let $e_N$ be the $N^{\text{th}}$ standard basis vector in $\mathbf{R}^N$, and define maps $\pi_1$ and $\pi_2$ from $G_N/\Gamma_N$ to $\mathcal{M}_N$, assigning $\pi_1(g\Gamma_N) = $ counting measure on $g\mathbf{Z}^N\setminus\{0\}$ and $\pi_2(g\Gamma_N) = $ counting measure on $g\Gamma_N e_N$. Let $\mu^j = \pi_j(\mu_N)$, where as before $\mu_N$ is normalized Haar measure on $G_N/\Gamma_N$. The Siegel relations (0.9)–(0.10) are tantamount to the statement that $\mu^1$ and $\mu^2$ are Siegel measures with $c(\mu^1) = 1$ and $c(\mu^2) = 1/\zeta(N)$.

To obtain Theorem 0.5 as a consequence of Theorem 0.12 it is only necessary to observe that by Masur's Theorem 0.1 the assignment to $x \in X$ of the counting measure $\nu_x$ on $\Pi(x)$ satisfies $\nu_x \in \mathcal{M}_2$. The fact that $\Pi(gx) = g\Pi(x)$, $x \in X$, $g \in G_2$ implies that $\nu_{gx} = g\nu_x$ and the image $\mu_0 \in \mathcal{P}(\mathcal{M}_2)$ of the measure $\mu$ is invariant and ergodic, i.e., a Siegel measure. The fact $\Pi(x) = -\Pi(x)$, by construction, implies $\nu_x \in \mathcal{M}_2^e$. Therefore, Parts I–III of Theorem 0.12 imply the corresponding parts of Theorem 0.5. A more complete discussion will be found in Sections 11–12.

As is illustrated by the Siegel theorem itself, a single ergodic action may give rise to more than one Siegel measure. This is especially true in the context of Theorem 0.5. With notations as in Theorem 0.5 define for each $x \in X$ and $0 \leq s < 1$ a set $\Pi(s,x) \subseteq \Pi(x)$ consisting of those vectors, with multiplicities, which arise from periodic cylinders of area $> s$. ($\Pi(0,x) = \Pi(x)$.) One finds $\Pi(s,gx) = g\Pi(s,x)$, $g \in G$, and by analogy with the preceding paragraph the map $x \to \nu_{x,s} = $ counting measure on $\Pi(s,x)$ determines a Siegel measure $\mu_s$ and a constant $c(\mu_s) \leq c(\mu_0)$ for which the conclusions of Theorem 0.5 remain true, i.e., with $\Pi(s,x)$ and $c(\mu_s)$ in place of $\Pi(x)$ and $c(\mu)$ respectively. The function $s \to c(\mu_s)$ is continuous from the right on $[0,1)$, but when $\mu$ is supported on an orbit, the function has finite range. If $(X, \mathcal{B}, \mu)$ is a component of a "stratum" of quadratic differentials equipped with its "Liouville measure" ([12], [21], [26], [22]), then with one trivial exception in genus one the function $s \to c(\mu_s)$ is continuous, positive and monotone decreasing, to zero, on $[0,1)$ (Theorem 13.3).

An important tool for the proof of Theorem 0.12 is Theorem 5.12, containing a basic identity which is derived in Sections 3–5. To describe this let $K = SO(N)$ and $A^+ = \{a = \text{diag}(a_1, \ldots, a_N) \mid \det a = 1, a_1 > \cdots > a_N > 0\}$. If $m_K(dk)$ is normalized Haar measure on $K$, and if $B = B(0,1)$ is the unit ball in $\mathbf{R}^N$, then

(0.16)
$$\int\limits_{\mathbf{R}^N}\int\limits_K \chi_B(akx)m_K(dk)\nu(dx)$$
$$= \frac{2\sigma_{N-1}}{\sigma_N}\int_0^1 \frac{N_\nu\left(\frac{\tau}{a_N}\right)}{\left(\frac{\tau}{a_N}\right)^N}(1-\tau^2)^{(N-3)/2}d\tau + O\left(\left(\frac{a_N}{a_{N-1}}\right)^{2/3}M(\nu)\right).$$

The identity (0.16) is used in Section 6 to prove the existence and finiteness of the constant $c(\mu)$ in Part I of Theorem 0.12, i.e., for (0.13). Given a Siegel measure $\mu$, the fact that the right side of (0.16) is bounded for each $\nu \in \mathcal{M}_N$ is combined with a corollary to a mean ergodic theorem, Theorem 2.6, to establish that $\hat{\chi}_B \in L^1(\mu)$. It is not difficult then to infer that $\hat{\psi} \in L^1(\mu)$ for each $\psi \in C_c(\mathbf{R}^N)$. Now (0.13) with $c(\mu) < \infty$ follows from uniqueness properties of Lebesgue measure.

A second application of (0.16), in an altered form, occurs in the proof of Part II of Theorem (0.12). One uses Part I, i.e., (0.13), to obtain

$$\int\limits_{\mathcal{M}_N} \frac{N_\nu(t)}{t^N}\mu(d\nu) = c(\mu)\frac{\sigma_N}{N} \qquad (0 < t < \infty).$$

This relation is used in Remark 5.21 to replace, at certain stages of the proof of Theorem 5.12, a pointwise error $M(\nu)$ by an $L^1(\mu)$ error $c(\mu)\frac{\sigma_N}{N}$. The result is an $L^1(\mu)$ error estimate $O\left(\left(\frac{a_N}{a_{N-1}}\right)^{2/3}c(\mu)\right)$ in (0.16) (Theorem 5.23). A version of the Wiener tauberian theorem is used then to establish (0.14) and Part II of Theorem 0.12 (Theorem 5.28).

The proof of Part III of Theorem 0.12 makes use of Part II and a result below which serves as a "Weyl criterion" for establishing that a net of even, locally finite Borel measures on $\mathbf{R}^N$ converges to Lebesgue measure in the $C_c(\mathbf{R}^N)$ topology. The Weyl criterion, Theorem 10.1, turns on the case $\nu_2 = $ Lebesgue measure of Theorem 9.4, here stated as

THEOREM 0.17. *Let $\nu_1, \nu_2 \in \mathcal{M}_N^e$ be such that $\nu_1(E) = \nu_2(E)$ for every ellipsoid $E$ centered at 0. Then $\nu_1 = \nu_2$.*

Theorems 0.3 and 0.17 will have as one corollary the following.

THEOREM 0.18. *Let $(X,q)$ and $c$ be as in Theorem 0.3. If $\psi \in C_c(\mathbf{R}^2)$, then*

(0.19) $$\lim_{R\to\infty}\frac{1}{R^2}\sum_{v\in\Pi(q)}\psi\left(\frac{v}{R}\right) = c\int\limits_{\mathbf{R}^2}\psi(u)du.$$

Let $q$ be such that $\Gamma(q)$ is a lattice, and define a zeta function $\zeta_q(s) = \sum_{v \in \Pi(q)} |v|^{-s}$, $\operatorname{Re} s > 2$. It is established in [23] that $\zeta_q(2s)$ is an entire-holomorphic linear combination of the Eisenstein series associated to the cusps of $\Gamma(q)$. From this one infers (a) $\zeta_q(\cdot)$ is entire meromorphic and (b) the constant $c$ which appears in (0.4) and (0.19) is essentially the residue of $\zeta_q(2s)$ at the simple pole $s = 1$. Explicit calculation of $c$ is possible in certain instances ([23], [24], [28]). In Section 15 we apply Theorems 5.19 and 10.1 and the equidistribution theorem of Eskin-McMullen ([3]) to prove Theorem 0.3 (and its consequence Theorem 0.18) without recourse to the theory of Eisenstein series. The role of [3] is to verify a property which we call "regularity" and which is motivated by [3], [5] and [17]. Briefly stated for the context of $(X, \mathcal{B}, \mu)$ in Theorem 0.5, a point $x \in X$ is $\mu$-regular if $\lim_{g \to \infty} g(m_K * \delta_x) = \mu$ in a suitable topology. When $x$ is regular, it develops that for all $0 < s < 1$ and $g \in G$ the set $g\Pi(s, x)$ has asymptotic growth $c(\mu_s)\pi R^2$ (Theorem 15.10). The prevalence of regularity in homogeneous space settings gives some hope for its genericity in the context of Theorem 0.5.

To further illustrate the "Eisenstein series free" approach to (0.4) and (0.19) we observe in Section 16, Theorem 16.1, that if $\Gamma = -\Gamma$ is a nonuniform lattice in $G = \operatorname{SL}(2, \mathbf{R})$, and if $v \in \mathbf{R}^2$ is such that $\Gamma v$ is a discrete set, then $\Gamma v$ satisfies (0.4) and (0.19) for a finite constant $c = c(\Gamma, v)$. Moreover, if $\Lambda$ is the isotropy group of $v$ in $\Gamma$, there is a number $t = t(\Lambda, v)$ such that $c(\Gamma, v) = 2t^2(\pi \operatorname{Vol} G/\Gamma)^{-1}$. One feature of the derivation of the formula for $c(\Gamma, v)$ is that it will not depend upon knowledge of meromorphic continuation of the Eisenstein series $E(z, s)$, $\operatorname{Im} z > 0$, $\operatorname{Re} s > 1$, which is associated to $(\Gamma, \Lambda)$. We shall give a direct proof that for each $z$ the function $(s - 1)E(z, s)$ has nontangential limit $(\operatorname{Vol} G/\Gamma)^{-1}$ at $s = 1$ from the half plane $\operatorname{Re} s > 1$. Of course, this implies the known fact that the residue of $E(z, \cdot)$ at $s = 1$ is $(\operatorname{Vol} G/\Gamma)^{-1}$ ([19], [7] (p. 224). The author thanks M. Wolf for providing the latter reference.)

Section 14 is devoted to the issue of pointwise a.e. convergence in Parts II–III of Theorem 0.12. If $N > 2$, and if we assume of the Siegel measure $\mu$ that $\hat{\psi} \in L^2(\mu)$ for all $\psi \in C_c(\mathbf{R}^N)$, then estimates in [9], Chapter V, are used to prove (0.14)–(0.15) are true for $\mu$-a.e. $\nu$. When $N = 2$, the same statement is true if one also assumes the representation of $G = \operatorname{SL}(2, \mathbf{R})$ on the orthocomplement of the constants does not almost have invariant vectors (cf. [9]).

Let $m = m(du)$ be Lebesgue measure on $\mathbf{R}^N$. If $c \geq 0$, define $S^c : \mathcal{M}_N \to \mathcal{M}_N$ by $S^c \nu = \nu + cm$. $S^c$ is equivariant relative to the action of $G_N$ on $\mathcal{M}_N$. In particular, if $\mu \in \mathcal{P}(\mathcal{M}_N)$ is a Siegel measure, then $S^c_* \mu$ is also a Siegel measure. We shall call a Siegel measure $\mu$ *singular* if $\nu \perp m$ for $\mu$-a.e. $\nu$. In the theorem to follow the point mass at the zero measure ($\nu \equiv 0$) is considered to be a singular Siegel measure.

THEOREM 0.20. *If $\mu$ is a Siegel measure, there exist $c \geq 0$ and a singular Siegel measure $\mu^s$ such that $\mu = S_*^c \mu^s$. In particular, if $\nu \prec m$ for $\mu$-a.e. $\nu$, then $\mu$ is a point mass at $cm$ for some $c \geq 0$.*

Theorem 0.20 is proved in Section 6 (Theorem 6.10). A second characterization of the point mass at $cm$ ($c \geq 0$) will be given in Section 8: Call $\nu$ *scale invariant* if $\nu(\lambda E) = |\lambda|^N \nu(E)$ for all Borel sets $E$ and real numbers $\lambda$. If $\mu$ is a Siegel measure such that $\nu$ is scale invariant for $\mu$-a.e. $\nu$, then $\mu$ is a point mass at $cm$ for some $c \geq 0$ (Theorem 8.6).

Work on this project was begun during a stay at the Laboratoire de Mathematiques Discretes with the kind support of Université d'Aix Marseille 2 (June, 1995). Indeed, the thought that Siegel's Theorem might be relevant, at least in spirit, to the study of periodic trajectories for quadratic differentials was provoked by a lecture on [4], at Luminy, by G.A. Margulis. The ideas in [4], [3] and [5] have been important to us.

The author wishes to thank M. Boshernitzan for useful conversations in connection with this work.

## 1. A mean ergodic theorem

Let $G$ be a semisimple analytic group with finite center and no compact factors. $G$ is a finite extension of a product of noncompact simple groups, $G \xrightarrow{\rho} (G_1 \times \cdots \times G_r)$, $\rho = (\rho_1, \ldots, \rho_r)$. The notation $g \to_s \infty$ is understood to mean $\rho_j(g) \to \infty$, $1 \leq j \leq r$.

If $\alpha = \{\alpha_n\}$ is a sequence in $G$, $U_\alpha$ shall denote the set of $g \in G$ such that the sequence $\{\alpha_n^{-1} g \alpha_n\}$ has the identity ($e$) for a cluster point. $N_\alpha$ is the least closed subgroup which contains $U_\alpha$. In the special case that $\alpha_n = b^n$, $n \geq 1$, we set $N_b = N_\alpha$ and $U_b = U_\alpha$, recalling that (1) $N_b = U_b$ and (2) for any sequence $\alpha$ there exists $b$ such that $N_\alpha = N_b$ (cf. [27]). Recall that $N_b$ is *totally unbounded* if $\rho_j(N_b) \neq \{e\}$, $1 \leq j \leq r$.

Let $V$ be a Banach space with norm $\|\cdot\|$, and let $\pi(\cdot)$ be a bounded, strongly continuous representation of $G$ on $V$. We make the standing assumption that there exist a bounded projection $P_\pi$ onto the subspace of invariant vectors and that $P_\pi \circ \pi(\cdot) = P_\pi$.

*Definition* 1.1. The representation $\pi$ above shall be called *admissible* if whenever $b \in G$ is such that $N_b$ is totally unbounded, then for each $v \in V$ the orbit $\pi(N_b)v$ has $P_\pi v$ in its norm closed convex hull.

Let $K$ be a fixed maximal compact subgroup of $G$, and let $m_K(dk)$ denote normalized Haar measure on $K$. $Q_\pi$ denotes the natural projection on the $K$-invariant vectors, defined by a Bochner integral as

$$(1.2) \qquad Q_\pi v = \int_K (\pi(k)v) m_K(dk).$$

With notations and definitions fixed above, we can state the mean ergodic theorem for admissible representations:

THEOREM 1.3. *Let $G$ be a semisimple analytic group with finite center and no compact factors. Let $\pi$ be an admissible representation of $G$ on a Banach space $V$, and let $Q_\pi$ be defined by (1.2) for a fixed choice of maximal compact subgroup $K$. Then*

$$(1.4) \qquad \lim_{g \to_s \infty} Q_\pi \pi(g) = P_\pi$$

*holds in the strong operator topology.*

*Proof.* Let $K$ be as in the statement of the theorem, and let $G = K A_c^+ K$ be a fixed Cartan decomposition. If $g = k_1 a k_2$, then $Q_\pi \pi(g) = Q_\pi \pi(a) \pi(k_2)$. The compactness of $K$ and boundedness and strong continuity of the representation $\pi$ combine to reduce (1.4) to

$$\lim_{\substack{a \to_s \infty \\ a \in A_c^+}} Q_\pi \pi(a) = P_\pi$$

in the strong operator topology. In fact, it is sufficient to prove

$$(1.4') \qquad \lim_{\substack{a \to_s \infty \\ a \in A_c^+}} \|Q_\pi \pi(a) v\| = 0 \qquad (v \in V, \; P_\pi v = 0).$$

Let $\theta(\cdot)$ denote the Cartan involution of $G$ which fixes $K$, and use the same notation for $\theta$ on the Lie algebra $\mathfrak{G}$. In order to establish (1.4') it is sufficient to consider sequences $\alpha = \{\alpha_n\} \subseteq A_c^+$ such that $\alpha_n \to_s \infty$ and $N_\alpha = \{h \in G \mid \lim_{n \to \infty} \alpha_n^{-1} h \alpha_n = e\}$. Then $\theta(N_\alpha) = N_{\alpha^{-1}}$ is totally unbounded and of the form $N_{\alpha^{-1}} = N_b$ for some $b \in G$. Letting $v \in V$ such that $P_\pi v = 0$ and $\epsilon > 0$ be given, we shall first apply the hypothesis of admissibility to find a probability measure $\xi$ on $N_{\alpha^{-1}}$ such that $\xi$ has compact support and

$$(1.5) \qquad \left\| \int_{N_{\alpha^{-1}}} \pi(h) v \xi(dh) \right\| < \epsilon.$$

Let $\mathfrak{n}_\alpha$, $\mathfrak{n}_{\alpha^{-1}}$ denote the Lie algebras of $N_\alpha$, $N_{\alpha^{-1}}$, respectively. $\log(\cdot)$ denotes the inverse to the exponential map where it is naturally defined. Set up a continuous map $k : N_{\alpha^{-1}} \to K$ as

$$k(h) = \exp\bigl(\log h + \theta(\log h)\bigr) \qquad (h \in N_{\alpha^{-1}}).$$

If we then set $k(h, n) = k(\alpha_n h \alpha_n^{-1})$, we have $\lim_{n \to \infty} k(h, n) = e$ locally uniformly on $N_{\alpha^{-1}}$. Define $g(h, n)$, implicitly, by

$$(1.6) \qquad hg(h, n) = \alpha_n^{-1} k(h, n) \alpha_n$$
$$= \exp(\log h + \alpha_n^{-2} \theta(\log h) \alpha_n^2).$$

Since $\theta(\log h) \in \mathfrak{n}_\alpha$, we also have $\lim_{n \to \infty} g(h, n) = e$ locally uniformly on $N_{\alpha^{-1}}$.

Let $\gamma_n$ be the image of the probability measure $\xi$ on $N_{\alpha^{-1}}$ under the map $h \to k(h, n)$. If $\|\pi\| = \sup_{g \in G} \|\pi(g)\|_{op}$, where $\| \cdot \|_{op}$ denotes operator norm, then $\|\pi\| < \infty$ by our boundedness assumption. If $h \in N_{\alpha^{-1}}$ and $v \in V$, define

$$(1.7) \qquad \delta(h, n, v) = \|\pi\| \left\| \Big(I - \pi(g(h, n))\Big) v \right\|.$$

Then $\lim_{n \to \infty} \delta(h, n, v) = 0$ locally uniformly in $h \in N_{\alpha^{-1}}$, $v$ fixed. We observe that if $h \in N_{\alpha^{-1}}$, then

$$(1.8) \qquad \|\pi(\alpha_n h) v - \pi(\alpha_n h g(h, n)) v\| \leq \delta(h, n, v).$$

Define $\delta(n, v)$ by

$$(1.9) \qquad \delta(n, v) = \sup_{h \in \text{sppt } \xi} \delta(h, n, v).$$

Since $\xi$ has compact support, $\lim_{n \to \infty} \delta(n, v) = 0$. Finally, by definition of $Q_\pi$ and (1.6)–(1.9) we have

$$\begin{aligned} Q_\pi \pi(\alpha_n) v &= \int_K \pi(k \alpha_n) v m_K(dk) \\ &= \int_K \int_K \pi(k k' \alpha_n) v \gamma_n(dk') m_K(dk) \\ &= \int_K \int_{N_{\alpha^{-1}}} \pi(k k(h, n) \alpha_n) v \xi(dh) m_K(dk) \\ &= \int_K \int_{N_{\alpha^{-1}}} \pi(k \alpha_n h g(h, n)) v \xi(dh) m_K(dk) \\ &= 0(\delta_n) + \int_K \int_{N_{\alpha^{-1}}} \pi(k \alpha_n h) v \xi(dh) m_K(dk) \\ &= 0(\delta_n) + 0(\|\pi\| \epsilon). \end{aligned}$$

Since $n$, then $\epsilon$ are arbitrary, we have $\lim_{n \to \infty} \|Q_\pi \pi(\alpha_n) v\| = 0$, and the theorem follows.

## 2. Applications of the mean ergodic theorem

We continue to suppose $G$ is a semisimple analytic group with finite center and no compact factors. Let $\mathcal{W} = \mathcal{W}(G)$ be the Banach algebra of continuous weakly almost periodic functions on $G$. Each $f \in \mathcal{W}$ is bounded and left and right uniformly continuous, and in particular the right regular representation is strongly continuous on $\mathcal{W}(G)(\pi(g)f(\cdot) = f(\cdot g))$. Denote by $\mathcal{E}(\cdot)$ the Eberlein mean on $\mathcal{W}$ ([2], [1], [18]). $\mathcal{E}(f)$ is the unique constant such that the set of (say) right translates of $f$ has the corresponding constant function in its (sup) norm convex hull. It is proved in [27] that

$$\lim_{g \to_s \infty} f(hg) = \mathcal{E}(f) \tag{2.1}$$

is true pointwise and, therefore, in the weak topology of $\mathcal{W}(G)$. It follows from the Banach-Mazur theorem that if $b \in G$ is such that $N_b$ is totally unbounded, then the constant function $\mathcal{E}(f)$ is in the norm closed convex hull of the orbit $\pi(N_b)f$. Therefore, the right regular representation on $\mathcal{W}(G)$ is admissible. From Theorem 1.3 we conclude

THEOREM 2.2. *Let $G$ be a semisimple analytic group with finite center and no compact factors. If $f \in \mathcal{W}(G)$, then for any maximal compact subgroup $K$*

$$\lim_{g \to_s \infty} \sup_{h \in G} \left| \int_K f(hkg) m_K(dk) - \mathcal{E}(f) \right| = 0. \tag{2.3}$$

Let $V$ be a reflexive Banach space, and let $\pi$ be a bounded strongly continuous representation of $G$ on $V$. The coefficients of $\pi$ belong to $\mathcal{W}(G)$ and if $f_{v_1, v_2^*}(g) = \langle \pi(g)v_1, v_2^* \rangle$, $v_1 \in V$, $v_2 \in V^*$, then $\mathcal{E}(f_{v_1, v_2^*}) = \langle P_\pi v_1, v_2^* \rangle$, where $P_\pi$ is a bounded equivariant projection on the subspace of invariant vectors. It now follows from (2.1) that if $b \in G$ is such that $N_b$ is totally unbounded, then for each $v \in V$ $P_\pi v$ is in the weak closure of the orbit $\pi(N_b)v$. The Banach-Mazur theorem implies $\pi$ is admissible, and we have

THEOREM 2.4. *Let $G$ be a semisimple analytic group with finite center and no compact factors. If $\pi$ is a bounded strongly continuous representation of $G$ on a reflexive Banach space $V$, and if $Q_\pi$ is defined in terms of a fixed maximal compact subgroup $K$, then*

$$\lim_{g \to_s \infty} Q_\pi \pi(g) = P_\pi \tag{2.5}$$

*in the strong operator topology.*

Our final application of Theorem 1.3 is in ergodic theory *per se*. Let $(X, \mathcal{B}, \mu)$ be a probability space, and let $G$ be represented there as a measurable group of measure preserving transformations. ("Measurable group" is understood to mean the pairing $(g, x) \to gx$ is measurable.) For our Banach space $V$ we take $\mathcal{L}^p(X, \mathcal{B}, \mu)$. Define $(\pi(g)f)(x) = f(g^{-1}x)$, $f \in \mathcal{L}^p$. If $\mathcal{B}_I$ is the $\sigma$-algebra of invariant measurable sets, i.e., $B \in \mathcal{B}_I$ when $\mu(gB \triangle B) = 0$, $g \in G$, then $P_\pi(\cdot) = E(\cdot \mid \mathcal{B}_I)$ is a contractive invariant projection.

THEOREM 2.6. *Let $G$ be a semisimple analytic group with finite center and no compact factors. Let $(X, \mathcal{B}, \mu)$ be a probability space, and let $G$ be there represented as a measurable group of measure preserving transformations. If $P_\pi = E(\cdot \mid \mathcal{B}_I)$ is the conditional expectation operator, then*

$$(2.7) \qquad \lim_{g \to_s \infty} \left\| \int_K f(gkx) m_K(dk) - E(f \mid \mathcal{B}_I) \right\|_p = 0$$

*for all $1 \le p < \infty$ and $f \in \mathcal{L}^p(X, \mathcal{B}, \mu)$.*

*Proof.* The integral in (2.7) is $(Q_\pi \pi(g^{-1})f)(x)$. If $1 < p < \infty$, apply Theorem 2.4. The *case* $p = 1$ follows by a standard argument from the case $p > 1$. □

Theorem 2.6 provides a useful criterion for integrability:

COROLLARY 2.8. *Let the notations and assumptions be as in Theorem 2.6. If $f \ge 0$ is measurable, and if for $\mu$-a.e. $x$*

$$(2.9) \qquad \limsup_{g \to_s \infty} \int_K f(gkx) m_K(dk) < \infty$$

*then $E(f \mid \mathcal{B}_I)(\cdot) < \infty$ a.e., where a.e. $x$*

$$(2.10) \qquad E(f \mid \mathcal{B}_I)(x) = \lim_{T \to \infty} E(\min(f, T) \mid \mathcal{B}_I)(x).$$

*In particular, if $\mu$ is ergodic for the $G$ action, (2.9) implies $f \in \mathcal{L}^1$.*

*Proof.* Define $f_T(x) = \min(f(x), T)$. Choose a sequence $g_n \to_s \infty$ such that $\lim_{n \to \infty} \int f_T(g_n kx) m_K(dk) = E(f_T \mid \mathcal{B}_I)(x)$ a.e. For $\mu$-a.e. $x$ we have

$$\infty > \limsup_{g \to_s \infty} \int f(gkx) m_K(dk)$$
$$\ge \limsup_{g \to_s \infty} \int f_T(gkx) m_K(dk)$$
$$\ge \lim_{n \to \infty} \int f_T(g_n kx) m_K(dk)$$
$$= E(f_T \mid \mathcal{B}_I)(x).$$

Now let $T \to \infty$ to obtain the desired result. □

*Remark* 2.11. For a certain class of groups $G$ Nevo [17] has obtained a mean ergodic theorem such as Theorem 2.6 for the family (in the present notation) $Q_\pi \pi(g) Q_\pi$. Nevo's purpose is for use in a much more delicate pointwise ergodic theorem.

## 3. Spherical integrals

This section is concerned with integrals similar to those in Section 2 of [4]. Our purpose is to establish first an elementary upper bound and then to motivate the calculation of derivatives in Section 4.

If $N > 1$, $\mathcal{A}_N^+$ shall denote the semigroup of diagonal matrices

$$\mathcal{A}_N^+ = \{\lambda = \mathrm{diag}(\lambda_1, \ldots, \lambda_N) \mid \lambda_1 > \lambda_2 > \cdots > \lambda_N > 0\}.$$

It is not required that $\det \lambda = 1$, $\lambda \in \mathcal{A}_N^+$.

We set $K = \mathrm{SO}(N)$ and let $m_K(\cdot)$ denote normalized Haar measure on $K$. If $S_{N-1}$ is the unit sphere in $\mathbf{R}^N$, and if $e_N$ is the $N^{\mathrm{th}}$ standard basis vector, then $S_{N-1} = K e_N$. If $\psi \geq 0$ is a Borel function on $\mathbf{R}^N$, define $\mathcal{C}\psi$ on $\mathcal{A}_N^+$ by

$$(3.1) \qquad \mathcal{C}\psi(\lambda) = \int_K \psi(\lambda k e_N) m_K(dk) \qquad (\lambda \in \mathcal{A}_N^+).$$

In view of the symmetry of $\lambda K e_N = \lambda S_{N-1}$ we shall always suppose $\psi$ is an *even* function. If $\sigma_N$ is the area of $S_{N-1}$, and if $S_{N-1}^+ = \{u \in S_{N-1} \mid u_N > 0\}$, then for an even Borel function $\psi \geq 0$ we have

$$(3.2) \qquad \mathcal{C}\psi(\lambda) = \frac{2}{\sigma_N} \int_{S_{N-1}^+} \psi(\lambda u) \frac{du_1 \wedge \cdots \wedge du_{N-1}}{u_N}.$$

In what follows $\psi$ is assumed to be, in addition to Borel and even, bounded with compact support. Fix $R_0 < \infty$ such that $\psi$ vanishes outside the ball ($B(0, R_0)$) of radius $R_0$ centered at $0$. If $\lambda \in \mathcal{A}_N^+$ is such that $\lambda_N \geq R_0$, then $\mathcal{C}\psi(\lambda) = 0$. Accordingly, we restrict attention to $\mathcal{A}_N^+(R_0) = \{\lambda \in \mathcal{A}_N^+ \mid 0 < \lambda_N \leq R_0\}$.

Next, fix $1 \leq k \leq N-1$, and define

$$\mathcal{A}_N^+(k, R_0) = \{\lambda \in \mathcal{A}_N^+(R_0) \mid \lambda_k \geq R_0 > \lambda_{k+1}\}.$$

Introduce a coordinate

$$(w, v) = (\lambda_1 u_1, \ldots, \lambda_k u_k, u_{k+1}, \ldots, u_N) = L_k u$$

on $S_{N-1}^+$, and observe that if $\lambda' = (\lambda_{k+1}, \ldots, \lambda_N)$

$$(3.3)$$
$$\lambda_1 \ldots \lambda_k \mathcal{C}\psi(\lambda) = \frac{2}{\sigma_N} \int_{L_k S_{N-1}^+} \frac{\psi(w, \lambda' v) dw_1 \wedge \cdots \wedge dw_k \wedge dv_{k+1} \wedge \cdots \wedge dv_{N-1}}{v_N}$$

where $(w, v) = L_k u$ as above. The right side of (3.3) is uniformly bounded on $\mathcal{A}_N^+(k, R_0)$. In particular, we have

LEMMA 3.4. *Let $\psi$ be an even, bounded Borel function with compact support, and let $\mathcal{C}\psi$ be defined on $\mathcal{A}_N^+$ by* (3.1). *There exists a constant $C =$*

$C(\psi) < \infty$ such that

$$(3.5) \qquad \mathcal{C}(\psi)(\lambda) \leq \frac{C(\psi)}{\lambda_1 \ldots \lambda_{N-1}} \qquad (\lambda \in \mathcal{A}_N^+).$$

In order to motivate a calculation in the next section assume of $\psi$ that for each $0 < \lambda_N < \infty$ almost all $z \in \mathbf{R}^{N-1}$ are such that $(z, \lambda_N)$ is a point of continuity of $\psi$. With this assumption we have for $H(\lambda) = \lambda_1 \ldots \lambda_{N-1} \mathcal{C}(\psi)(\lambda)$

$$\lim_{\substack{\lambda' \in \mathcal{A}_N^+ \\ \lambda'_k \to \infty, k < N \\ \lambda'_N \to \lambda_N}} H(\lambda') = \frac{2}{\sigma_N} \mathcal{R}\psi(\lambda_N)$$

where $\mathcal{R}\psi$ is the restricted Radon transform

$$\mathcal{R}\psi(\lambda_N) = \int_{\mathbf{R}^{N-1}} \psi(w, \lambda_N) dw_1 \wedge \cdots \wedge dw_{N-1}.$$

This is evident from (3.3) with $k = N$ since for large $(\lambda'_1, \ldots, \lambda'_{N-1})$, $v = u_N$ must be close to 1 and $\lambda'_N v$ close to $\lambda_N$.

If $\psi$ is of class $C^1$, and if we set $\psi_t(u) = \psi(tu)$, then because $\nabla \psi(y) \cdot y = \frac{d}{dt}\psi(ty)\big|_{t=1}$, one has

$$\mathcal{R}(\nabla\psi(x) \cdot x)(\lambda_N) = (1-N)\mathcal{R}\psi(\lambda_N) + (\mathcal{R}\psi)'(\lambda_N)\lambda_N.$$

From (3.2) we have

$$\nabla_\lambda \mathcal{C}\psi(\lambda) \cdot \lambda = \frac{2}{\sigma_N} \int_{S_{N-1}^+} \nabla\psi(\lambda u) \cdot \lambda u \frac{du_1 \wedge \cdots \wedge du_{N-1}}{u_N}$$

and then

$$(3.6) \quad \lim_{\substack{\lambda'_k \to \infty, k < N \\ \lambda'_N \to \lambda_N}} \lambda'_1 \ldots \lambda'_{N-1} \left(\nabla_{\lambda'}\mathcal{C}\psi(\lambda') \cdot \lambda'\right) = \frac{2}{\sigma_N} \mathcal{R}(\nabla\psi(x) \cdot x)(\lambda_N)$$

$$= \frac{2}{\sigma_N}\left[(1-N)\mathcal{R}\psi(\lambda_N) + (\mathcal{R}\psi)'(\lambda_N) \cdot \lambda_N\right].$$

Set aside the question of differentiability and let $\psi = \chi_B$, where $B = B(0,1)$ is the unit ball in $\mathbf{R}^N$. By direct calculation $\mathcal{R}\chi_B(\lambda_N) = \frac{\sigma_{N-1}}{N-1}(1-\lambda_N^2)^{(N-1)/2}$ and $(\mathcal{R}\chi_B)'(\lambda_N) = -\sigma_{N-1}(1-\lambda_N^2)^{(N-3)/2}\lambda_N$ when $0 \leq \lambda_N < 1$. Then since

$$\frac{2}{\sigma_N}\left[(1-N)\frac{\sigma_{N-1}}{N-1}(1-\lambda_N^2)^{(N-1)/2} - \sigma_{N-1}(1-\lambda_N^2)^{(N-3)/2}\lambda_N^2\right]$$

$$= \frac{2\sigma_{N-1}}{\sigma_N}(1-\lambda_N^2)^{(N-3)/2},$$

it is a formal consequence of (3.6) that

$$(3.7) \quad \lim_{\substack{\lambda'_k \to \infty, k < \infty \\ \lambda'_N \to \lambda_N \\ \lambda' \in \mathcal{A}_N^+}} \lambda'_1 \ldots \lambda'_{N-1} \nabla \mathcal{C}\chi_B(\lambda') \cdot \lambda' = -\frac{2\sigma_{N-1}}{\sigma_N}(1-\lambda_N^2)^{(N-3)/2},$$

which holds for $0 \leq \lambda_N < 1$. We shall verify a version of (3.7), with bounds, in Section 4.

## 4. Calculations for the ball

Let $1 < c < \infty$, and redefine
$$\mathcal{A}_N^+(c) = \{\lambda \in \mathcal{A}_N^+ \mid \lambda_{N-1} > c > 1 > \lambda_N\}.$$
To fix notations in this section we define $F_N(\lambda) = \mathcal{C}\chi_B(\lambda)$, where $B = B(0,1)$. In this section we shall prove

THEOREM 4.1. *If $1 < c < \infty$, there exists a constant $\eta(c) < \infty$ such that if $\lambda \in \mathcal{A}_N^+(c)$, then*

(4.2) $$\left| \lambda_1 \ldots \lambda_{N-1}(\nabla F_N(\lambda) \cdot \lambda) + \frac{2\sigma_{N-1}}{\sigma_N}(1 - \lambda_N^2)^{(N-3)/2} \right|$$
$$\leq \frac{\eta(c)}{\lambda_{N-1}^2} \frac{2\sigma_{N-1}}{\sigma_N}(1 - \lambda_N^2)^{(N-3)/2}.$$

*As the notation suggests $F_N(\cdot)$ is continuously differentiable on $\mathcal{A}_N^+(c)$.*

The inductive proof involves some calculations, the first being used to establish (4.2) when $N = 2$.

Let $\lambda \in \mathcal{A}_2^+(c)$ for some $c > 1$, and define $a(\lambda) = \left(\frac{1-\lambda_2^2}{\lambda_1^2 - \lambda_2^2}\right)^{1/2}$. Declare $\sigma_1 = 2$ so that the constant $\frac{2\sigma_1}{\sigma_2} = \frac{2}{\pi}$. Now
$$F_2(\lambda) = \frac{1}{\pi} \int_{-a(\lambda)}^{a(\lambda)} \frac{du}{(1-u^2)^{1/2}} = \frac{2}{\pi} \sin^{-1} a(\lambda).$$

Then by a short calculation
$$\lambda_1 \nabla F_2(\lambda) \cdot \lambda = -\frac{2}{\pi} \frac{\lambda_1}{(\lambda_1^2 - 1)^{1/2}(1 - \lambda_2^2)^{1/2}}$$

and one finds

(4.3) $$\eta(c) = 2c^2 \left(\frac{c}{(c^2-1)^{1/2}} - 1\right).$$

Now suppose $N > 2$ and the theorem has been established for $N - 1$. Fix $1 < c < \infty$, and assume below that $\lambda \in \mathcal{A}_N^+(c)$, i.e., that $\lambda_{N-1} > c > 1 > \lambda_N$.

We introduce some coordinates and other quantities. Assume $|u_1| < 1/\lambda_1$ ($< 1$), $u \in S_{N-1}^+$

(4.4) $$v_j(\lambda, u) = \frac{u_j}{(1-u_1^2)^{1/2}}$$
$$(2 \leq j \leq N)$$
$$\mu_j(\lambda, u) = \frac{(1-u_1^2)^{1/2}}{(1-\lambda_1^2 u_1^2)^{1/2}} \lambda_j.$$

While $\{v_j\}$ depend upon all coordinates of $u$, one should note that $\{\mu_j\}$ involve only the first, i.e., $u_1$. Define $a(\lambda)$, by analogy with $(a(\lambda), N = 2)$ above, as

$$a(\lambda) = \left(\frac{1 - \lambda_N^2}{\lambda_1^2 - \lambda_N^2}\right)^{1/2}.$$

We have

(4.5)
$$a(\lambda) < 1/\lambda_1$$
$$\mu_N(\lambda, u_1) \geq 1, \quad |u_1| \geq a(\lambda)$$
$$\mu_j(\lambda, u_1) > \lambda_j, \quad 2 \leq j \leq N, \quad |u_1| < 1/\lambda_1.$$

By abuse of notation we view $\mu = (\mu_2, \ldots, \mu_N)$ as an element of $\mathbf{R}^{N-1}$. By (4.5) we have $\mu \in \mathcal{A}_{N-1}^+(c)$ when $\lambda \in \mathcal{A}_N^+(c)$ and $|u_1| < a(\lambda)$. When $a(\lambda) \leq |u_1| < 1/\lambda_1$, we observe $F_{N-1}(\mu) = 0$.

By definition $v(\lambda, u) \in S_{N-2}^+$ when $u \in S_{N-1}^+$, and $\sum_{j=2}^{N} \mu_j^2 v_j^2 < 1$ precisely when $\sum_{j=1}^{N} \lambda_j^2 u_j^2 < 1$. Since (4.4) implies

$$\frac{du_2 \wedge \cdots \wedge du_{N-1}}{u_N} = \left(1 - u_1^2\right)^{\frac{N-1}{2}} \frac{dv_2 \wedge \cdots \wedge dv_{N-1}}{v_N}$$

we have

$$F_N(\lambda) = \frac{2}{\sigma_N} \int_{-a(\lambda)}^{a(\lambda)} \left(1 - u_1^2\right)^{\frac{N-1}{2}} \int_{S_{N-2}^+} \chi_{B_{N-1}}(\mu v) \frac{dv_2 \wedge \cdots \wedge dv_{N-1}}{v_N} du_1$$
$$= \frac{\sigma_{N-1}}{\sigma_N} \int_{-a(\lambda)}^{a(\lambda)} \left(1 - u_1^2\right)^{\frac{N-1}{2}} F_{N-1}(\mu(\lambda, u_1)) du_1.$$

Since $F_{N-1}(\mu(\lambda, \pm a(\lambda))) = 0$, it is possible to differentiate under the integral to find

$$\nabla F_N(\lambda) \cdot \lambda = \frac{\sigma_{N-1}}{\sigma_N} \int_{-a(\lambda)}^{a(\lambda)} \left(1 - u_1^2\right)^{\frac{N-1}{2}} \nabla_\lambda \left(F_{N-1} \circ \mu(\lambda, u_1)\right) \cdot \lambda du_1.$$

Use the relations $\frac{\partial \mu_i}{\partial \lambda_j} = \frac{\delta_{ij}}{\lambda_i}$, $2 \leq i, j \leq N$, and $\frac{\partial \mu_j}{\partial \lambda_1} = \frac{-\lambda_1 u_1^2}{1 - \lambda_1^2 u_1^2} \mu_j$, $2 \leq j \leq N$, to establish

$$\nabla_\lambda F_{N-1} \circ \mu(\lambda, u_1) \cdot \lambda = \frac{1}{1 - \lambda_1^2 u_1^2} \nabla_\mu F_{N-1}(\mu) \cdot \mu.$$

We now have

(4.6) $$\nabla F_N(\lambda) \cdot \lambda = \frac{\sigma_{N-1}}{\sigma_N} \int_{-a(\lambda)}^{a(\lambda)} \frac{(1 - u_1^2)^{\frac{N-1}{2}}}{(1 - \lambda_1^2 u_1^2)} \nabla_\mu F_{N-1}(\mu) \cdot \mu \, du_1.$$

Definition (4.3) implies

(4.7) $$\lambda_1 \ldots \lambda_{N-1} = \lambda_1 \frac{(1 - \lambda_1^2 u_1^2)^{\frac{N-2}{2}}}{(1 - u_1^2)^{\frac{N-2}{2}}} \mu_2 \ldots \mu_{N-1}.$$

Multiply the two sides of (4.6) by $\lambda_1 \ldots \lambda_{N-1}$ and use (4.7) to find
(4.8)
$$\lambda_1 \ldots \lambda_{N-1} \nabla F_N(\lambda) \cdot \lambda$$
$$= \lambda_1 \frac{\sigma_{N-1}}{\sigma_N} \int_{-a(\lambda)}^{a(\lambda)} (1-u_1^2)^{1/2}(1-\lambda_1^2 u_1^2)^{\frac{N-4}{2}} \mu_2 \ldots \mu_{N-1} \nabla F_{N-1}(\mu) \cdot \mu \, du_1.$$

As we have observed, $\mu(\lambda, u_1) \in \mathcal{A}_{N-1}^+(c)$ and $\mu_j > \lambda_j$ in the range of $u_1$ of interest. We apply the induction hypothesis to write

$$(4.9) \quad \mu_2 \ldots \mu_{N-1} \nabla F_{N-1}(\mu) \cdot \mu = -\frac{2\sigma_{N-2}}{\sigma_{N-1}}(1-\mu_N^2)^{\frac{N-4}{2}}\left(1 + \frac{E(\mu)}{\mu_{N-1}^2}\right)$$

with $|E(\mu)| \le \eta_{N-1}(c)$. Using $\mu_{N-1} > \lambda_{N-1}$, the second summand on the right side in (4.9) satisfies

$$(4.10) \quad \left|\frac{E(\mu)}{\mu_{N-1}^2}\right| \le \frac{\eta_{N-1}(c)}{\lambda_{N-1}^2}.$$

We must now compute the principal part of the integral (4.8). To this end observe that

$$(1-\mu_N^2)^{\frac{N-4}{2}} = \frac{(1-\lambda_N^2)^{\frac{N-4}{2}}}{(1-\lambda_1^2 u_1^2)^{\frac{N-4}{2}}}\left(1-\left(\frac{u_1}{a(\lambda)}\right)^2\right)^{\frac{N-4}{2}}.$$

The principal part of (4.8), i.e., without the factor $1 + \frac{E(\mu)}{\mu_{N-1}^2}$ in (4.9), is now

$$(4.11) \quad -(1-\lambda_N^2)^{\frac{N-4}{2}}\lambda_1 \frac{2\sigma_{N-2}}{\sigma_N} \int_{-a(\lambda)}^{a(\lambda)}(1-u_1^2)^{1/2}\left(1-\left(\frac{u_1}{a(\lambda)}\right)^2\right)^{\frac{N-4}{2}} du_1.$$

Substitute $u_1 = a(\lambda)\tau$ so that (4.11) becomes

$$(4.12) \quad (1-\lambda_N^2)^{\frac{N-3}{2}} \frac{\lambda_1}{(\lambda_1^2 - \lambda_N^2)^{1/2}}\left(\frac{-2\sigma_{N-2}}{\sigma_N}\right)\int_{-1}^{1}(1-a^2(\lambda)\tau)^{1/2}(1-\tau^2)^{\frac{N-4}{2}} d\tau.$$

It is clear that $\frac{\lambda_1}{(\lambda_1^2 - \lambda_N^2)^{1/2}} = 1 + O\left(\frac{1}{\lambda_{N-1}^2}\right)$ and $(1-a^2(\lambda)\tau^2)^{1/2} = 1 + O\left(\frac{1}{\lambda_{N-1}^2}\right)$ hold uniformly for $\lambda \in \mathcal{A}_N^+(c)$ and $-1 \le \tau \le 1$. Make the substitution $\tau^2 = t$ in (3.12) to find
(4.13)
$$(1-\lambda_N^2)^{\frac{N-3}{2}}\left(\frac{-2\sigma_{N-2}}{\sigma_N}\right)\int_{-1}^{1}(1-\tau^2)^{\frac{N-4}{2}} d\tau\left(1 + O\left(\frac{1}{\lambda_{N-1}^2}\right)\right)$$
$$= (1-\lambda_N^2)^{\frac{N-3}{2}}\left(\frac{-2\sigma_{N-2}}{\sigma_N}\right) 2\int_0^1 (1-t)^{\frac{N-2}{2}-1} t^{\frac{1}{2}-1} \frac{dt}{2}\left(1 + O\left(\frac{1}{\lambda_{N-1}^2}\right)\right)$$
$$= (1-\lambda_N^2)^{\frac{N-3}{2}}\left(\frac{-2\sigma_{N-2}}{\sigma_N}\right) B\left(\frac{N-2}{2}, \frac{1}{2}\right)\left(1 + O\left(\frac{1}{\lambda_{N-1}^2}\right)\right)$$

where $B(\cdot,\cdot)$ is Euler's beta function. One knows $\sigma_{N-1} = \sigma_{N-2} B\left(\frac{N-2}{2}, \frac{1}{2}\right)$. Therefore, if we combine (4.9) and (4.13)

$$\lambda_1 \ldots \lambda_{N-1}(\nabla F_N(\lambda) \cdot \lambda) = -\frac{2\sigma_{N-1}}{\sigma_N}(1 - \lambda_N^2)^{\frac{N-3}{2}}\left(1 + O\left(\frac{1}{\lambda_{N-1}^2}\right)\right)$$

is true uniformly on $\mathcal{A}_N^+(c)$, $c > 1$. This implies the statement of Theorem 4.1.

## 5. A basic identity

Let $\nu$ be a Borel measure on $\mathbf{R}^N$. If $B(0, R)$ is the open ball of radius $R$ about 0, define the growth function

(5.1) $$N_\nu(R) = \nu(B(0, R)).$$

Then define

(5.2) $$M(\nu) = \sup_{0 < R < \infty} \frac{N_\nu(R)}{R^N}.$$

We shall deal with measures under the assumption $M(\nu) < \infty$. For certain formulae it is required only that $\frac{N_\nu(R)}{R^N}$ be bounded on every interval $(0, T)$, $T < \infty$.

Let $A_N^+ = \mathcal{A}_N^+ \cap G$, $G = \mathrm{SL}(N, \mathbf{R})$. If $a \in A_N^+$ and $R > 0$, define $Ra = \mathrm{diag}(Ra_1, \ldots, Ra_N) \in \mathcal{A}_N^+$. For later reference we record for $\lambda = \lambda(R, a)$, $\lambda_j = Ra_j$,

(5.3) $$\lambda_1 \ldots \lambda_{N-1} = \frac{R^{N-1}}{a_N}.$$

Let $F_N(\lambda)$ be as defined in Section 4. If $M(\nu) < \infty$, and if $B$ is the unit ball in $\mathbf{R}^N$, then

(5.4) $$\int_{\mathbf{R}^N} \int_K \chi_B(akx) m_K(dk) \nu(dx) = \int_0^\infty F_N(Ra) dN_\nu(R).$$

This relation holds because the inner integral on the left side is a radial function of $x$ (namely $F_N(\|x\|a)$).

LEMMA 5.5. *Assume of $a \in A_N^+$ that*

(5.6) $$\frac{a_N}{a_{N-1}} < \frac{1}{2}.$$

*Then*

(5.7) $$\int_0^{\frac{2}{a_{N-1}}} F_N(Ra) dN_\nu(R) = O\left(M(\nu)\frac{a_N}{a_{N-1}}\right)$$

*with $O(\cdot)$ uniform in $\nu$ and $a \in A_N^+$ satisfying* (5.6).

*Proof.* Set $C(\chi_B) = C_N$ in (3.5). From (3.5) and (5.3) we obtain $F_N(Ra) \leq C_N \frac{a_N}{R^{N-1}}$. Replace $F_N(Ra)$ by $C_N \frac{a_N}{R^{N-1}}$ on the left side in (5.7) and integrate by parts:

$$\int_0^{\frac{2}{a_{N-1}}} F_N(Ra) dN_\nu(R) \leq C_N \int_0^{\frac{2}{a_{N-1}}} \frac{a_N}{R^{N-1}} dN_\nu(R)$$

$$= C_N \frac{N_\nu(R)}{R^N}(a_N R) \Big|_0^{\frac{2}{a_{N-1}}} + C_N(N-1) \int_0^{\frac{2}{a_{N-1}}} \frac{a_N N_\nu(R)}{R^N} dR$$

$$= O\left(M(\nu) \frac{a_N}{a_{N-1}}\right)$$

as claimed.

Next, we take up the integral on the right-hand side of (5.4) over the interval $\left(\frac{2}{a_{N-1}}, \infty\right)$. If $Ra_N \geq 1$, then $F_N(Ra) = 0$, and therefore we study the integral over the finite interval $\left(\frac{2}{a_{N-1}}, \frac{1}{a_N}\right)$. By (5.6) this interval is not empty. Integrate by parts, taking into account the fact $F_N(a_N^{-1} a) = 0$ and the estimate which resulted in (5.7):

(5.8)
$$\int_{\frac{2}{a_{N-1}}}^{\frac{1}{a_N}} F_N(Ra) dN_\nu(R) = O\left(M(\nu) \frac{a_N}{a_{N-1}}\right) - \int_{\frac{2}{a_{N-1}}}^{\frac{1}{a_N}} N_\nu(R) \frac{d}{dR} F_N(Ra) dR.$$

Now by Theorem 4.1 and (5.3)

$$\frac{d}{dR} F_N(Ra) = \frac{1}{R} \nabla F_N(Ra) \cdot (Ra)$$

$$= -\frac{a_N}{R^N} \frac{2\sigma_{N-1}}{\sigma_N} (1 - R^2 a_N^2)^{\frac{N-3}{2}} \left(1 + O\left(\frac{1}{R^2 a_{N-1}^2}\right)\right)$$

with $O\left(\frac{1}{R^2 a_{N-1}^2}\right)$ uniform for $Ra_{N-1} \geq 2$ (say), as is true on the interval of integration in (5.8).

The second summand on the right side of (5.8) can be rewritten as

$$(5.9) \quad a_N \left(\frac{2\sigma_{N-1}}{\sigma_N}\right) \int_{\frac{2}{a_{N-1}}}^{\frac{1}{a_N}} \frac{N_\nu(R)}{R^N} (1 - R^2 a_N^2)^{\frac{N-3}{2}} \left(1 + O\left(\frac{1}{R^2 a_{N-1}^2}\right)\right) dR.$$

Substitute $\tau = Ra_N$ in (5.9) to obtain the integral

$$(5.10) \quad \frac{2\sigma_{N-1}}{\sigma_N} \int_{\frac{2a_N}{a_{N-1}}}^{1} \frac{N_\nu\left(\frac{\tau}{a_N}\right)}{\left(\frac{\tau}{a_N}\right)^N} (1-\tau^2)^{\frac{N-3}{2}} \left(1 + O\left(\frac{1}{\tau^2}\left(\frac{a_N}{a_{N-1}}\right)^2\right)\right) d\tau.$$

To deal with the $O(\cdot)$ term divide the interval of integration in (5.10) into intervals $\left(\frac{2a_N}{a_{N-1}}, \left(\frac{2a_N}{a_{N-1}}\right)^{2/3}\right)$ and $\left(\left(\frac{2a_N}{a_{N-1}}\right)^{2/3}, 1\right)$. The $O(\cdot)$ term is at most $1/4$ on the first interval. If we now assume

$$(5.11) \quad \frac{a_N}{a_{N-1}} < \frac{1}{4}$$

the integral over $\left(\frac{2a_N}{a_{N-1}}, \left(\frac{2a_N}{a_{N-1}}\right)^{2/3}\right)$ is $O\left(M(\nu)\left(\frac{a_N}{a_{N-1}}\right)^{2/3}\right)$ uniform in $\nu$ and $a \in A_N^+$ satisfying (5.11). As for the second interval the $O(\cdot)$ term under the integral is $O\left(\left(\frac{a_{N-1}}{a_N}\right)^{4/3}\left(\frac{a_N}{a_{N-1}}\right)^2\right) = O\left(\left(\frac{a_N}{a_{N-1}}\right)^{2/3}\right)$, and the $O(\cdot)$ term again contributes $O\left(\left(\frac{a_N}{a_{N-1}}\right)^{2/3} M(\nu)\right)$ to the integral. Finally, with $a \in A_N^+$ constrained by (5.11) we have

$$\frac{2\sigma_{N-1}}{\sigma_N} \int_{\left(\frac{2a_N}{a_{N-1}}\right)^{2/3}}^{1} \frac{N_\nu\left(\frac{\tau}{a_N}\right)}{\left(\frac{\tau}{a_N}\right)^N} (1-\tau^2)^{\frac{N-3}{2}} d\tau$$
$$= O\left(M(\nu)\left(\frac{a_N}{a_{N-1}}\right)^{2/3}\right) + \frac{2\sigma_{N-1}}{\sigma_N} \int_0^1 \frac{N_\nu\left(\frac{\tau}{a_N}\right)}{\left(\frac{\tau}{a_N}\right)^N} (1-\tau^2)^{\frac{N-3}{2}} d\tau.$$

By collecting results we have

THEOREM 5.12. *If $\nu$ is a Borel measure on $\mathbf{R}^N$ such that $M(\nu) < \infty$, then for $a \in A_N^+$ constrained by (5.11) there is the uniform estimate*

$$(5.13) \quad \int_{\mathbf{R}^N} \int_K \chi_B(akx) m_K(dk) \nu(dx)$$
$$= \frac{2\sigma_{N-1}}{\sigma_N} \int_0^1 \frac{N_\nu\left(\frac{\tau}{a_N}\right)}{\left(\frac{\tau}{a_N}\right)^N} (1-\tau^2)^{\frac{N-3}{2}} d\tau + O\left(M(\nu)\left(\frac{a_N}{a_{N-1}}\right)^{2/3}\right).$$

Define $h_N$ on $\mathbf{R}^+$ by $h_N(\tau) = \tau(1-\tau^2)^{\frac{N-3}{2}} \chi_{(0,1)}(\tau)$. The factor $\tau$ guarantees $h_N \in L^1\left(\mathbf{R}^+, \frac{d\tau}{\tau}\right)$. The Fourier transform, $\hat{h}_N$, is given by

$$\hat{h}_N(c) = \int_0^\infty h_N(\tau) \tau^{ic} \frac{d\tau}{\tau}$$

$$= \int_0^1 (1-\tau^2)^{\frac{N-3}{2}} \tau^{ic} d\tau$$

$$= \frac{1}{2} \int_0^1 (1-s)^{\frac{N-3}{2}} s^{ic/2 - 1/2} ds$$

$$= \frac{1}{2} B\left(\frac{N-1}{2}, \frac{1+ic}{2}\right)$$

$$= \frac{1}{2} \frac{\Gamma\left(\frac{N-1}{2}\right) \Gamma\left(\frac{1+ic}{2}\right)}{\Gamma\left(\frac{N+ic}{2}\right)}.$$

Since $N > 0$, $\hat{h}_N(c) \neq 0$, $c \in \mathbf{R}$. Since $\sigma_N = \sigma_{N-1} B\left(\frac{N-1}{2}, \frac{1}{2}\right)$, the function $g_N = \frac{2\sigma_{N-1}}{\sigma_N} h_N$ is integrable with $\hat{g}_N(0) = 1$. Define $\psi_\nu \in L^\infty\left(\mathbf{R}^+, \frac{d\tau}{\tau}\right)$ by $\psi_\nu(t) = t^N N_\nu\left(\frac{1}{t}\right)$, $t > 0$. The integral on the right-hand side of (5.13) is $g_N * \psi_\nu(a_N)$. Since $\hat{g}_N(c) \neq 0$, $c \in \mathbf{R}$, the Wiener tauberian theorem implies that if $\lim_{a_N \to 0} g_N * \psi_\nu(a_N) = \ell = \hat{g}_N(0)\ell$, then for every $g \in L^1\left(\mathbf{R}^+, \frac{d\tau}{\tau}\right)$ $\lim_{a_N \to 0} g * \psi_\nu(a_N) = \hat{g}(0)\ell$. In particular, if we set $g(t) = t\chi_{(0,1)}(t)$ and use Theorem 5.12, we obtain

THEOREM 5.14.  *Assume $M(\nu) < \infty$. If*

(5.15) $$\lim_{\substack{a \to \infty \\ a \in A_N^+}} \int_{\mathbf{R}^N} \int_K \chi_B(akx) m_K(dk) \nu(dx) = \ell$$

*then*

(5.16) $$\lim_{T \to \infty} \frac{1}{T} \int_0^T \frac{N_\nu(R)}{R^N} dR = \ell.$$

LEMMA 5.17.  *Let $\varphi > 0$ be defined on $\mathbf{R}^+$, and assume there exists $\alpha > 0$ such that $\varphi(u)u^\alpha$ is monotone nondecreasing. For all $t > 0$ and $\lambda > 1$ we have*

(5.18) $$\varphi^{-1}(0, \lambda\varphi(t)) \supseteq \left[\frac{t}{\lambda^{1/\alpha}}, t\right]$$

$$\varphi^{-1}\left(\frac{\varphi(t)}{\lambda}, \infty\right) \supseteq \left[t, \lambda^{1/\alpha} t\right].$$

*Proof.* If $t/\lambda^{1/\alpha} \leq s \leq t$, the assumption on $\varphi$ implies $\varphi(t)t^\alpha \geq \varphi(s)s^\alpha \geq \frac{t^\alpha}{\lambda}\varphi(s)$, and therefore $\varphi(s) \leq \lambda\varphi(t)$. The first line of (5.18) is true and the second line follows by a similar argument. □

We shall apply Theorem 5.14 and Lemma 5.17 to establish

THEOREM 5.19. *If $M(\nu) < \infty$, and if (5.15) is true, then*

$$\lim_{R \to \infty} \frac{N_\nu(R)}{R^N} = \ell. \tag{5.20}$$

*Proof.* Let $\varphi(R) = \frac{N_\nu(R)}{R^N}$ and $\alpha = N$ in Lemma 5.17. Set $A_r^s = \frac{1}{s-r}\int_r^s \varphi(t)dt$, $0 \leq r < s < \infty$. Note that $\ell < \infty$ by Theorem 5.12. Let $\ell^- = \liminf_{R \to \infty} \varphi(R)$. If $\ell^- < \ell$, set $a = \frac{\ell + \ell^-}{2} < \ell$ and fix $\lambda > 1$ such that $\lambda a < \ell$. The lemma implies that if $\varphi(t) < a$, then $\varphi < \lambda a < \ell$ on $\left[\frac{1}{\lambda^{1/\alpha}}t, t\right]$. Let $s = \frac{1}{\lambda^{1/\alpha}}t$. Then $A_0^t = \frac{s}{t}A_0^s + \left(1 - \frac{s}{t}\right)A_s^t \leq \frac{s}{t}A_0^s + \left(1 - \frac{s}{t}\right)\lambda a$. If $t \to \infty$ in such a way that $\varphi(t) < a$, it follows that $\ell \leq \frac{1}{\lambda^{1/\alpha}}\ell + \left(1 - \frac{1}{\lambda^{1/\alpha}}\right)\left(\frac{\ell + \ell^-}{2}\right) < \ell$, a contradiction. We conclude $\ell^- \geq \ell$. By an analogous argument $\limsup_{R \to \infty} \varphi(R) \leq \ell$, and the theorem follows. □

*Remark* 5.21. In most of the applications the measure $\nu$ will itself be a point in a probability space $(\mathcal{M}_N, \mathcal{B}_N, \mu)$ (§6). There will exist a constant $c(\mu) < \infty$ such that

$$\int_{\mathcal{M}_N} \frac{N_\nu(t)}{t^N} \mu(d\nu) = c(\mu)\frac{\sigma_N}{N}, \quad 0 < t < \infty. \tag{5.22}$$

The analysis which lead to (5.13) (Theorem 5.12) was conducted for a single $\nu$, $M(\nu) < \infty$, and involved replacing expressions of the form $\frac{N_\nu(t)}{t^N}$, in certain places by $M(\nu)$. If the issue in (5.13) is an error estimate for an equation in $L^1(\mu)$, one may integrate over $\mathcal{M}_N$ in the same places, applying the Fubini theorem where necessary, and thus replace $\frac{N_\nu(t)}{t^N}$ by its $L^1(\mu)$ norm ((5.22)), i.e, by $c(\mu)\frac{\sigma_N}{N}$. With this modification Theorem 5.12 may be restated as

THEOREM 5.23. *Let $(\mathcal{M}_N, \mathcal{B}_N, \mu)$ be as in Remark 5.21. If $a \in A_N^+$ is constrained by (5.11), there is the estimate*

$$\left\| \int_{\mathbf{R}^N} \int_K \chi_B(akx) m_K(dk) \nu(dx) - \frac{2\sigma_{N-1}}{\sigma_N} \int_0^1 \frac{N_\nu\left(\frac{\tau}{a_N}\right)}{\left(\frac{\tau}{a_N}\right)^N}(1-\tau^2)^{\frac{N-3}{2}}d\tau \right\|_1$$
$$= O\left(c(\mu)\left(\frac{a_N}{a_{N-1}}\right)^{2/3}\right). \tag{5.24}$$

Theorem 5.23 is a restatement of Theorem 5.12 in the context of Remark 5.21, i.e., for the measure spaces $(\mathcal{M}_N, \mathcal{B}_N, \mu)$ which satisfy (5.22). We shall now give corresponding replacements for Theorems 5.14 and 5.19.

THEOREM 5.25. *Let $(\mathcal{M}_N, \mathcal{B}_N, \mu)$ be as in Remark* 5.21. *If $\ell \in \mathbf{R}^+$ is such that*

$$\text{(5.26)} \quad \lim_{\substack{a \to \infty \\ a \in A_N^+}} \left\| \int_{\mathbf{R}^N} \int_K \chi_B(akx) m_K(dk) \nu(dx) - \ell \right\|_{1,\mu} = 0$$

*then*

$$\text{(5.27)} \quad \lim_{T \to \infty} \left\| \frac{1}{T} \int_0^T \frac{N_\nu(R)}{R^N} dR - \ell \right\|_{1,\mu} = 0.$$

*Proof.* As in the proof of Theorem 5.14 the integrand in the second term in (5.24) is expressed as $\psi_\nu * g_N(a_N)$. But now $\psi_\nu(t)$ is viewed as a function on $\mathbf{R}^+$ with values in $L^1(\mathcal{M}_N, \mathcal{B}_N, \mu)$. By (5.22) we have for each $t$ the relation $\|\psi_\nu(t)\|_{1,\mu} = c(\mu)\frac{\sigma_N}{N}$. That is, $\psi_\nu(t)$ is a *bounded* function from $\mathbf{R}^+$ to $L^1(\mathcal{M}_N, \mathcal{B}_N, \mu)$. It follows that the set $\{g \in L^1(\mathbf{R}^+, \frac{d\tau}{\tau}) \mid \lim_{a_N \to 0} \|\psi_\nu * g(a_N) - \ell \hat{g}(0)\|_{1,\mu} = 0\}$ is a closed ideal. Since $\hat{g}_N$ is never 0, this ideal is all of $L^1(\mathbf{R}^+, \frac{d\tau}{\tau})$. As in the proof of Theorem 5.14 the choice $g(t) = t\chi_{(0,1)}(t)$ yields (5.27). The theorem is proved. □

Finally, we shall replace Theorem 5.19 in the context of Remark 5.21.

THEOREM 5.28. *Let $(\mathcal{M}_N, \mathcal{B}_N, \mu)$ be as in Remark* 5.21. *If* (5.26) *is true, then*

$$\text{(5.29)} \quad \lim_{R \to \infty} \left\| \frac{N_\nu(R)}{R^N} - \ell \right\|_{1,\mu} = 0.$$

*Of course, by* (5.22), $\ell = c(\mu)\frac{\sigma_N}{N}$.

*Proof.* By the proof of Theorem 5.19, the relations (5.18) from Lemma 5.17 and the fact (5.27) is also true in measure imply that (5.29) is true in measure. It is therefore sufficient to establish that the family $\left\{\frac{N_\nu(R)}{R^N} \mid R \gg 0\right\}$ is uniformly integrable. To this end let $c > 0$ and let $E \in \mathcal{B}_N$, $R > 0$ be such that

$$\text{(5.30)} \quad \int_E \frac{N_\nu(R)}{R^N} \mu(d\nu) \geq c.$$

If $R \leq S \leq 2R$, then $\frac{N_\nu(S)}{S^N} \geq \frac{1}{2^N}\frac{N_\nu(R)}{R^N}$, and therefore by the Fubini theorem and (5.30)

$$\text{(5.31)} \quad \int_E \frac{1}{2R} \int_0^{2R} \frac{N_\nu(s)}{s^N} ds \mu(d\nu) \geq \frac{c}{2^{N+1}}.$$

Since $\left\{\frac{1}{T}\int_0^T \frac{N_\nu(s)}{s^N} ds \mid T \gg 0\right\}$ is uniformly integrable, (5.30)–(5.31) imply that for every $\epsilon > 0$ there is a $\delta > 0$ such that if $R > 1$ and $\mu(E) < \delta$, then $c < \epsilon$ in (5.30). The theorem is proved. □

## 6. Siegel measures

If $N > 1$, we define $\mathcal{M}_N$ to be the set of Borel measures on $\mathbf{R}^N$ such that $M(\nu) < \infty$. If $\psi(\cdot)$ is a compactly supported bounded Borel function on $\mathbf{R}^N$, define $\hat{\psi}$ on $\mathcal{M}_N$ by duality,

$$\hat{\psi}(\nu) = \int_{\mathbf{R}^N} \psi(x) \nu(dx). \tag{6.1}$$

Endow $\mathcal{M}_N$ with the smallest topology such that $\hat{\psi} \in C(\mathcal{M}_N)$ when $\psi \in C_c(\mathbf{R}^N)$ (= continuous, compactly supported functions on $\mathbf{R}^N$). The following fact implies $\mathcal{M}_N$ is a countable union of compact metrizable spaces:

LEMMA 6.2. *If $c < \infty$ and $\mathcal{M}_N(c) = \{\nu \mid M(\nu) \leq c\}$, then $\mathcal{M}_N(c)$ is compact and metrizable. In particular, $\mathcal{M}_N$ is a standard Borel space* ([30]).

*Proof.* The elementary proof is left to the reader.

*Remark.* Of course, $\mathcal{M}_N$ is neither locally compact nor metrizable.

Let $G = \mathrm{SL}(N, \mathbf{R})$. If $A \in G$, $A$ determines a linear transformation of $\mathbf{R}^N$ which, as a continuous map, maps measures to measures. Since $A^{-1}\nu(E) = \nu(AE)$, $\nu \in \mathcal{M}_N$, $E$ Borel, $A \in G$, we have

$$\begin{aligned} N_{A^{-1}\nu}(R) &\leq N_\nu(\|A\|R) \\ &\leq \|A\|^N R^N M(\nu). \end{aligned}$$

That is, $M(A^{-1}\nu) \leq \|A\|^N M(\nu)$, and $G$ acts naturally upon $\mathcal{M}_N$. It is clear that $(A, \nu) \to A^{-1}\nu$ is continuous as a map from $G \times \mathcal{M}_N$ to $\mathcal{M}_N$.

Denote by $\mathcal{P}(\mathcal{M}_N)$ the set of Borel probability measures on $\mathcal{M}_N$. We introduce

*Definition* 6.3. An element $\mu \in \mathcal{P}(\mathcal{M}_N)$ is a *Siegel measure* if (a) $G\mu = \mu$ relative to the action $(A, \nu) \to A^{-1}\nu$ and (b) $\mu$ is ergodic relative to this action.

THEOREM 6.4. *Let $\mu$ be a Siegel measure. If $\psi \in C_c(\mathbf{R}^N)$, then $\hat{\psi} \in L^1(\mu)$.*

*Proof.* It is a consequence of Theorem 5.12 that if $0 \leq \psi < c\chi_B$, $B = B(0,1)$, $c < \infty$, and if we express $g \in G$ as $g = k_1 a k_2$, $a \in A_N^+$, $k_1, k_2 \in K = \mathrm{SO}(N)$, then (note $g \to \infty$ is the same as $g \to_s \infty$)

$$\limsup_{g \to \infty} \int_K \hat{\psi}(gk\nu) m_K(dk) \leq \limsup_{k_1 a k_2 = g \to \infty} c \int_{\mathbf{R}^N} \int_K \chi_B(k_1 a k_2 k x) m_K(dk) \nu(dx)$$

$$= \limsup_{\substack{a \to \infty \\ a \in A_N^+}} c \int_{\mathbf{R}^N} \int_K \chi_B(akx) m_K(dx) \nu(dx) \leq c C_N M(\nu) < \infty.$$

By Corollary 2.8 $\hat{\psi} \in L^1(\mu)$. For the general case define $T_\lambda \nu(E) = \frac{\nu(\lambda E)}{\lambda^N}$, $\lambda > 0$, where $\lambda E$ denotes the homothety by $\lambda$ on $\mathbf{R}^N$. Since $T_\lambda G = G T_\lambda$, the image $T_\lambda \mu$ of a Siegel measure under this map of $\mathcal{M}_N$ is again a Siegel measure. Since $\hat{\psi}(T_\lambda \nu) = \hat{\psi}_\lambda(\nu)$, where $\psi_\lambda(y) = \lambda^{-N}\psi(\lambda^{-1}y)$, and since $\psi_\lambda$ is supported on $B(0,1)$ when $\psi$ is supported on $B(0,1/\lambda)$, the first part of the argument implies $\hat{\psi} \in L^1(T_\lambda \mu)$ when $\psi$ is supported on $B(0,1/\lambda)$. But since $\mu = T_\lambda T_{\lambda^{-1}}\mu$ and $T_{\lambda^{-1}}\mu$ is Siegel, we have $\hat{\psi} \in L^1(\mu)$. The lemma is proved. □

THEOREM 6.5.  *If $\mu \in \mathcal{P}(\mathcal{M}_N)$ is a Siegel measure, there exists $c(\mu) < \infty$ such that for any Borel function $\psi \in L^1(\mathbf{R}^N, dx)$,*

$$\tag{6.6} \int_{\mathcal{M}_N} \hat{\psi}(\nu)\mu(d\nu) = c(\mu) \int_{\mathbf{R}^N} \psi(x)dx$$

*where $dx$ is Lebesgue measure.*

*Proof.* If $\psi \in C_c(\mathbf{R}^N)$, then $\hat{\psi} \in L^1(\mu)$. Define a functional $\Phi$ by

$$\tag{6.7} \Phi(\psi) = \int_{\mathcal{M}_N} \hat{\psi}(\nu)\mu(d\nu) \qquad (\psi \in C_c(\mathbf{R}^N)).$$

If $\psi^A(x) = \psi(A^{-1}x)$ and $\hat{\psi}^A(\nu) = \hat{\psi}(A^{-1}\nu)$, then $\widehat{\psi^A} = \hat{\psi}^A$. Since $\mu$ is invariant, we have $\Phi(\psi^A) = \Phi(\psi)$, $A \in G$, $\psi \in C_c(\mathbf{R}^N)$. Also, $\psi \geq 0$ implies $\Phi(\psi) \geq 0$. It follows then that there exist $a, b \geq 0$ such that

$$\Phi(\psi) = a\psi(0) + b \int_{\mathbf{R}^N} \psi(x)dx \qquad (\psi \in C_c(\mathbf{R}^N)).$$

Choose $\psi_k(x) = \chi_B(x)(1-\|x\|^2)^k$, $k > 0$, so that $0 \leq \psi_k \leq \chi_B$ and $\psi_k(x) \to \delta_{0x}$ pointwise. By definition of $M(\nu)$ we have $\hat{\psi}_k(\nu) \to 0$ for all $\nu$, and by the dominated convergence theorem

$$a = \lim \Phi(\psi_k) = 0.$$

It follows that we may take $b = c(\mu)$ so that (6.6) is true when $\psi \in C_c(\mathbf{R}^N)$. The extension to integrable Borel functions $\psi$ on $\mathbf{R}^N$ is straightforward and will be omitted (e.g., the set of $\psi$ which satisfy (6.6) is closed under monotone (integrable) limits). □

COROLLARY 6.8.  *Let $\mu \in \mathcal{P}(\mathcal{M}_N)$ be a Siegel measure. Then (5.22) is true. That is*

$$\tag{6.9} \int_{\mathcal{M}_N} \frac{N_\nu(t)}{t^N}\mu(d\nu) = c(\mu)\frac{\sigma_N}{N} \qquad (0 < t < \infty).$$

We conclude this section with a rudimentary structure theorem for Siegel measures. Let $m = m_N$ denote Lebesgue measure on $\mathbf{R}^N$. If $c \geq 0$, define $S^c : \mathcal{M}_N \to \mathcal{M}_N$ by

$$S^c \nu = \nu + cm \qquad (\nu \in \mathcal{M}_N,\ c \geq 0).$$

Since $GS^c = S^c G$, the induced map $S^c_*$, on $\mathcal{P}(\mathcal{M}_N)$ sends Siegel measures to Siegel measures. We call a Siegel measure $\mu$ *singular* if $\nu \perp m$ for $\mu$-a.e. $\nu$. If $c \geq 0$, $\eta_c$ denotes the point mass at $cm$. Trivially, $\eta_c \in \mathcal{P}(\mathcal{M}_N)$ is a Siegel measure. By Lemma 6.2 both $\mathcal{M}_N$ and $\mathcal{M}_N \times \mathbf{R}^N$ are standard Borel spaces. We now restate Theorem 0.20 of the introduction in the form

THEOREM 6.10. *If $\mu \in \mathcal{P}(\mathcal{M}_N)$ is a Siegel measure, then either* (a) $\mu = \eta_{c(\mu)}$ *or* (b) *there exist $c \geq 0$ and a singular Siegel measure $\mu^s$ such that $\mu = S^c_* \mu^s$.*

*Proof.* Identify $(\mathbf{R}^N \setminus \{0\}, m(du))$ with $(G/H_0, m_{G/H_0}(dgH_0))$ where $H_0$ is the isotropy group of an arbitrary but fixed vector $v_0 \in \mathbf{R}^N \setminus \{0\}$. If $\mu$ is a Siegel measure, then Moore's ergodicity theorem ([16]) implies $\mu$ is ergodic for the induced action of $H_0$. Then, according to Zimmer ([31, Theorem 4.2]; see also [30, Proposition 2.22]), $\mu \times m_{G/H_0}$ is ergodic for the $G$ action on $\mathcal{M}_N \times G/H_0$. Associate to $\mu$ a Borel measure $\lambda_\mu$ on $\mathcal{M}_N \times G/H_0$, defined by disintegration (cf. [6]) as

$$（6.11） \quad \lambda_\mu(F) = \int_{\mathcal{M}_N} \nu(F[\nu]) \mu(d\nu).$$

In (6.11) $F \subseteq \mathcal{M}_N \times G/H_0$ is a Borel set and $F[\nu] = \{gH_0 \mid (\nu, gH_0) \in F\}$, $\nu \in \mathcal{M}_N$. $\lambda_\mu$ is $\sigma$-finite (e.g., $\lambda_\mu(\mathcal{M}_N \times B(0,R)) = c(\mu) R^N \frac{\sigma_N}{N}$ by Theorem 6.5). Let $\lambda_\mu = \lambda_\mu^a + \lambda_\mu^s$ be the Lebesgue decomposition of $\lambda_\mu$ relative to $\mu \times m$, i.e., $\lambda_\mu^a \prec \mu \times m$ and $\lambda_\mu^s \perp \mu \times m$. Both $\lambda_\mu^a$ and $\lambda_\mu^s$ are $G$-invariant, and therefore the ergodicity of $\mu \times m$ implies $\lambda_\mu^a = c \cdot (\mu \times m)$ for some $c \geq 0$.

Since $\lambda_\mu^s \perp \mu \times m$, there exists a Borel set $E \subseteq \mathcal{M}_N \times G/H_0$ such that $\mu \times m(E) = 0 = \lambda_\mu^s(E^c)$. If $\nu \in \mathcal{M}_N$, define a Borel measure $\xi_\nu$ on $G/H_0$ by $\xi_\nu(A) = \nu(A \cap E[\nu])$. Since $\xi_\nu \leq \nu$ and $G/H_0 \cong \mathbf{R}^N \setminus \{0\}$, there is a natural sense in which $\xi_\nu \in \mathcal{M}_N$.

If $F \subseteq \mathcal{M}_N \times G/H_0$ is a Borel set, then by (6.11) and the choice of $E$, we have

$$（6.12） \quad \begin{aligned} \lambda_\mu^s(F) &= \lambda_\mu(E \cap F) \\ &= \int_{\mathcal{M}_N} \nu\big((E \cap F)[\nu]\big) \mu(d\nu) \\ &= \int_{\mathcal{M}_N} \xi_\nu(F[\nu]) \mu(d\nu). \end{aligned}$$

If $B \subseteq G/H_0$ is a Borel set, and if we set $F = \mathcal{M}_N \times B$ in (6.12), we find that $\nu \to \xi_\nu(B)$ is a Borel function on $\mathcal{M}_N$. From this we conclude that $\nu \to \hat{\psi}(\xi_\nu)$ is a Borel function on $\mathcal{M}_N$ for each $\psi \in C_c(\mathbf{R}^N)$, and therefore $R(\nu) = \xi_\nu$ is a Borel map of $\mathcal{M}_N$ to itself. In particular, (6.12) is a disintegration of $\lambda_\mu^s$ over $\mu$.

If $g \in G$, the fact $g\lambda_\mu^s = \lambda_\mu^s$ combines with the a.e. $\mu$ uniqueness of the representation (6.12) to imply $g\xi_\nu = \xi_{g\nu}$, $\mu$-a.e. $\nu$. We conclude that $R : \mathcal{M}_N \to \mathcal{M}_N$ is Borel and a.e. $\mu$ an equivariant map. In particular, $\mu^s = R_*\mu$ is a Siegel measure. By construction $\mu^s$ is singular. Since $\nu = cm + \xi_\nu$, a.e. $\nu$, we have $S^c R = \mathrm{Id}$ a.e. $\mu$ and $\mu = S_*^c \mu^s$. The theorem is proved. □

## 7. Asymptotic growth

If $\mathcal{B}_N$ is the Borel $\sigma$-algebra of $\mathcal{M}_N$, then Corollary 6.8 implies that for every Siegel measure $\mu$ the triple $(\mathcal{M}_N, \mathcal{B}_N, \mu)$ satisfies the hypothesis (5.22) of Remark 5.21. Collecting results from Theorems 2.6, 5.25 and 5.28 we have

THEOREM 7.1. *Let $N > 1$, and assume $\mu \in \mathcal{P}(\mathcal{M}_N)$ is a Siegel measure. Then*

$$\text{(7.2)} \qquad \lim_{R \to \infty} \left\| \frac{N_\nu(R)}{R^N} - c(\mu)\frac{\sigma_N}{N} \right\|_1 = 0.$$

*Remark* 7.3. If $N > 2$, and if the Siegel measure $\mu$ is such that $\hat{\psi} \in L^2(\mu)$ for all $\psi \in C_c(\mathbf{R}^N)$, then (7.2) also holds pointwise for $\mu$-a.e. $\nu$ (Theorem 14.11).

## 8. Special Siegel measures

Let $\nu \in \mathcal{M}_N$ be such that $T_\lambda \nu = \nu$, $\lambda > 0$, where $T_\lambda \nu(E) = \frac{\nu(\lambda E)}{\lambda^N}$ for Borel sets $E$. Such a $\nu$ has a unique expression in polar coordinates as

$$\text{(8.1)} \qquad \nu = \gamma_\nu \times (R^{N-1} dR)$$

for a finite Borel measure $\gamma_\nu$ on $S_{N-1}$. Since $N_\nu(R) = \gamma_\nu(S_{N-1})\frac{R^N}{N}$, we have

$$\text{(8.2)} \qquad \frac{N_\nu(R)}{R^N} = \frac{\gamma_\nu(S_{N-1})}{N}.$$

PROPOSITION 8.3. *Let $\mu$ be a Siegel measure such that $T_\lambda \nu = \nu$, $\lambda > 0$, for $\mu$-a.e. $\nu$. Define $h(\nu) = \gamma_\nu(S_{N-1})$, where $\gamma_\nu$ is defined by (8.1). Then if $c(\mu)$ is as in Theorem 6.5,*

$$\text{(8.4)} \qquad h(\nu) = c(\mu)\sigma_N \qquad (\mu - \text{a.e. } \nu).$$

*Proof.* Immediate from (7.2) and Theorem 7.1. □

LEMMA 8.5. *If $\nu \in \mathcal{M}_N$ is such that $T_\lambda \nu = \nu$, $\lambda > 0$, then for each $A \in G = \mathrm{SL}(N, \mathbf{R})$ $A^{-1}\nu$ has the same property. If $h(\nu) = \gamma_\nu(S_{N-1})$, then*

$$\text{(8.6)} \qquad h(A^{-1}\nu) = \int_{S_{N-1}} \frac{\gamma_\nu(dx)}{\|A^{-1}x\|^N}.$$

*Proof.* If $E \subseteq S_{N-1}$ is a Borel set, and if $\widehat{E} = \{tx \mid x \in E, 0 < t < 1\}$, then $\nu(\widehat{E}) = \frac{\gamma_\nu(E)}{N}$. If $E = S_{N-1}$, then $h(A^{-1}\nu) = \gamma_{A^{-1}\nu}(S_{N-1}) = N\nu(AB)$, $B =$ unit ball $= \widehat{S}_{N-1}$. If $x \in S_{N-1}$, then $x = \|A^{-1}x\| \left(A \frac{A^{-1}x}{\|A^{-1}x\|}\right)$, and therefore $AB$ contains the interval $\{tx \mid 0 < t < \frac{1}{\|A^{-1}x\|}\}$. It follows then that

$$h(A^{-1}\nu) = N\nu(AB)$$
$$= \int_{S_{N-1}} \frac{\gamma_\nu(dx)}{\|A^{-1}x\|^N}$$

as claimed. □

In the next theorem we shall assume $\mu$ is a Siegel measure on $\mathcal{M}_N^e$, the set of *even* elements of $\mathcal{M}_N$.

THEOREM 8.7. *Let $\mathcal{M}_N^e = \{\nu \in \mathcal{M}_N \mid \nu(-E) = \nu(E), E \text{ Borel}\}$. If $\mu$ is a Siegel measure on $\mathcal{M}_N^e$, and if $T_\lambda \nu = \nu$, $\lambda > 0$, for $\mu$-a.e. $\nu \in \mathcal{M}_N^e$, then $\mu$ is a point mass at $\nu = c(\mu)dx$, i.e.,*

(8.8) $$\mu = \eta_{c(\mu)}$$

*in the notation of Theorem* 6.10.

*Proof.* Let $\mathbf{P}_{N-1} = S_{N-1}/\pm 1$. Since $\mu$-a.e. $\nu$ is even, the measure $\gamma_\nu$ in (8.2) is even. Proposition 8.4 and Lemma 8.5 imply

(8.9) $$\int_{S_{N-1}} \frac{\gamma_\nu(dx)}{\|g^{-1}x\|^N} = c(\mu)\sigma_N \qquad (g \in G, \mu\text{-a.e. } \nu).$$

If $\sigma(dx)$ is Euclidean measure on $S_{N-1}$, the measure $c(\mu)\sigma(dx)$ is also even and satisfies (8.9). It is only necessary to establish that an even measure $\gamma_\nu$ on $S_{N-1}$ is uniquely determined by the integrals (8.9). Now if $[x] = \pm x$ is an element of the real projective space $\mathbf{P}_{N-1}$, the function $P(gK, [x]) = \|g^{-1}x\|^N$ is a Poisson kernel on $G/K \times \mathbf{P}_{N-1}$ for the Laplace-Beltrami operator on $G/K$ and the nonmaximal boundary $\mathbf{P}_{N-1}$. It is well known and easily proved that the Poisson integrals $\varphi(gK) = \int_{\mathbf{P}_{N-1}} P(gK, [x])\lambda(d[x])$, $\lambda$ a finite measure on $\mathbf{P}_{N-1}$, uniquely determine $\lambda$. (See Remark 8.10.) Since $\gamma_\nu$ and $c(\mu)\sigma(dx)$ are even measures on $S_{N-1}$ with the same Poisson integrals, their projections on $\mathbf{P}_{N-1}$ are equal. That is, $\gamma_\nu = c(\mu)\sigma(dx)$, as claimed.

*Remark* 8.10. Let $G = \mathrm{SL}(N, \mathbf{R})$, $K = \mathrm{SO}(N)$, and let $M$ be the group of diagonal elements of $K$. Then $B(G) = K/M$ is the Furstenberg boundary of $G$. The Poisson kernel on $G/K \times B(G)$ is

(8.11) $$P_o(gK, kM) = e^{-2\rho(H(g^{-1}k))}$$
$$= \prod_{j=1}^{N-1} \|g^{-1}k(e_1 \wedge \cdots \wedge e_j)\|^{-2}$$

where $e_1, \ldots, e_N$ is the standard basis for $\mathbf{R}^N$ (cf. [25, Section 4]). According to [11, Propositions 2.6–2.6′], a finite positive measure $\lambda$ on $B(G)$ is uniquely determined by its Poisson integral $P\lambda(gK) = \int_{B(G)} P_o(gK, kM)\lambda(dkM)$. Map $K/M$ to $\mathbf{P}_{N-1}$ by $kM \to [\pm ke_1]$. Let $K_{N-1} = \mathrm{SO}(N-1)$ be embedded in $\mathrm{SO}(N)$ by $u \to k(u) = \begin{pmatrix} 1 & 0 \\ 0 & u \end{pmatrix}$, and let $du$ be normalized Haar measure on $K_{N-1}$. The key relation is

$$(8.12) \qquad \frac{1}{\|g^{-1}ke_1\|^N} = \int_{K_{N-1}} P_o(gK, kk(u)M)du.$$

To establish this relation directly, let $g^{-1}k = k_1 an$ be expressed as a Iwasawa decomposition. In view of (8.11) we may suppose $k_1 = I$. Let $a = \mathrm{diag}(a_1, \ldots, a_N)$, and let $\hat{a}$, $\hat{n}$ be the $(N-1) \times (N-1)$ matrix consisting of the rows, columns $2, \ldots, N$ of $a$, $n$, respectively. We have $\det \hat{a} = a_1^{-1}$ and $\det \hat{n} = 1$. For all $j$ and $u \in K_{N-1}$

$$(8.13) \qquad g^{-1}kk(u)e_1 \wedge \cdots \wedge e_j = (a_1 e_1) \wedge \hat{a}\hat{n}u(e_2 \wedge \cdots \wedge e_j).$$

Let $\tilde{a} = a_1^{1/N-1} \hat{a}$ so that $\det \tilde{a} = 1$, and observe that (8.13) implies

$$(8.14) \qquad \prod_{j=1}^{N-1} \|g^{-1}kk(u)(e_1 \wedge \cdots \wedge e_j)\|^{-2}$$

$$= a_1^{-2(N-1)} a_1^{\sum_{j=1}^{N-1} \frac{2(j-1)}{N-1}} \prod_{j=2}^{N-1} \|\tilde{a}\hat{n}u(e_2 \wedge \cdots \wedge e_j)\|^{-2}$$

$$= a_1^{-N} \prod_{j=2}^{N-1} \|\tilde{a}\hat{n}u(e_2 \wedge \cdots \wedge e_j)\|^{-2}.$$

If $\tilde{g}^{-1} = \tilde{a}\hat{n} \in \mathrm{SL}(N-1, \mathbf{R})$, then since $a_1 = \|g^{-1}ke_1\|$, the last term in (8.13) is $\|g^{-1}ke_1\|^{-N} P'_o(\tilde{g}K_{N-1}, uM_{N-1})$. Here $P'_o$ is the Poisson kernel for $\mathrm{SL}(N-1, \mathbf{R})/K_{N-1} \times B(\mathrm{SL}(N-1, \mathbf{R}))$. Since the Poisson kernel has integral 1, (8.12) follows. The uniqueness statement in the proof of Theorem 8.7 is now a direct consequence of [11]. Simply lift the image of $\gamma_\nu$ on $\mathbf{P}_{N-1}$ to $B(G)$ with the help of $du$ and the map $kM \to [\pm ke_1]$.

## 9. A characterization of Lebesgue measure

This section is devoted to proof of the following theorem and its corollary below:

THEOREM 9.1. *Let $N > 1$, and let $\nu$ be a Borel measure on $\mathbf{R}^N$ such that (a) $\nu$ is even, i.e., $\nu(-U) = \nu(U)$ for every Borel set $U$ and (b) $\nu(E) = m(E)$* *for every ellipsoid $E$ with center $0$, where $m(\cdot)$ is Lebesgue measure. Then $\nu = m$.*

COROLLARY 9.2. *Let $\psi \in C_c(\mathbf{R}^N)$ and $\epsilon > 0$ be given. There exist ellipsoids $E_1, \ldots, E_r$ centered at $0$ and $\delta > 0$ such that if $\nu$ is an even Borel measure, and if $|\nu(E_j) - m(E_j)| < \delta$, then*

$$\left| \int \psi(x)\nu(dx) - \int \psi(x)m(dx) \right| < \epsilon.$$

*Proof of corollary.* Suppose the statement is false. There exist $\psi \in C_c(\mathbf{R}^N)$ and $\epsilon > 0$ such that for every finite set $\mathcal{E}$ of ellipsoids with center zero there is a Borel measure $\nu = \nu_{\mathcal{E}}$ with the properties

(9.3) $$|\nu(E) - m(E)| < 1/(\text{Card } \mathcal{E}) \quad (E \in \mathcal{E})$$
$$\left| \int \psi(x)\nu(dx) - \int \psi(x)m(dx) \right| \geq \epsilon.$$

Let $\mathcal{F}$ be the set of such finite sets $\mathcal{E}$, ordered by inclusion. The net $\{\nu_{\mathcal{E}}\}_{\mathcal{E} \in \mathcal{F}}$ is locally bounded, and therefore there is a subnet $\{\nu_{\mathcal{E}'}\}$ such that $\lim_{\mathcal{E}'} \nu_{\mathcal{E}'} = \nu$ exists in the $C_c(\mathbf{R}^N)$ topology. Use $rE$ to denote homothety of an ellipsoid $E$ by $r > 0$. For any $E$ and $r_1 < 1 < r_2$, we have $\{r_1 E, r_2 E\} \subseteq \mathcal{E}'$, large $\mathcal{E}'$. Now

$$\left| \nu_{\mathcal{E}'}(r_j E) - r_j^N m(E) \right| < 1/(\text{Card } \mathcal{E}').$$

This implies $r_1^N m(E) \leq \nu(E) \leq r_2^N m(E)$ for all $r_1 < 1 < r_2$, and therefore $\nu(E) = m(E)$ for every ellipsoid $E$ with center zero. Theorem 9.1 implies $\nu = m$ and, in particular,

$$\lim_{\mathcal{E}'} \int \psi(x)\nu_{\mathcal{E}'}(dx) = \int \psi(x)m(dx).$$

This contradicts (9.3), and the corollary is proved. □

We shall prove Theorem 9.1 in two steps. First, we shall assume

(9.4) $$\nu(dx) = \psi(x)m(dx)$$

where $\psi(\cdot)$ is uniformly bounded and continuous on $\mathbf{R}^N \setminus \{0\}$ and $\psi(-x) \equiv \psi(x)$. An approximation (convolution) argument is then used to reduce to the first case.

If $0 < \theta < \pi/2$, and if $\varphi_{N-1}(\theta)$ is the surface area of the set of $x \in S_{N-1}$ whose spherical distance from $e_1 = (1, 0, \ldots, 0)$ is less than $\theta$, then

(9.5) $$\varphi_{N-1}(\theta) = \frac{\sigma_{N-1}}{N-1} (\sin \theta)^{N-1}(1 + o(1))$$

where $o(1)$ is as $\theta \to 0$. (Recall that $\sigma_{N-1}$ is the surface area of $S_{N-2}$.)

We focus on the ray $\mathbf{R}^+ e_1$. If $0 < a < 1$ and $R > 0$, define $E(a, R)$ to be the ellipsoid

$$E(a, R) = \left\{(x, y) \mid x \in \mathbf{R}, \ y \in \mathbf{R}^{N-1}, \ x^2 + \frac{\|y\|^2}{a^2} < R^2\right\}.$$

For later reference we record

(9.6) $$m(E(a, R)) = a^{N-1} R^N \frac{\sigma_N}{N} = \nu(E(a, R)).$$

We shall study $\nu$ in polar coordinates

(9.7) $$\nu(d(R, x)) = \psi(Rx) R^{N-1} dR A(dx)$$

where $A(dx)$ is the Euclidean surface area measure on $S_{N-1}$. It is our goal to prove $\psi \equiv 1$.

If $0 < t < 1$, the intersection of $E(a, R)$ with the sphere $S(tR) = \{(x, y) \mid x^2 + \|y\|^2 = t^2 R^2\}$ has one or two components, one if $t \leq a$ and two if $t > a$. The contribution to $\nu(E(a, R))$ from values $t \leq a$ is $O(a^N)$. As we shall be letting $a \to 0$, it will be no loss to assume $a < t < 1$.

If $0 < x = tR \cos\theta$ and $\|y\| = tR \sin\theta$, then $(x, y) \in S(tR) \cap E(a, R)$ if and only if

$$\sin\theta < \frac{a}{t}\left(\frac{1 - t^2}{1 - a^2}\right)^{1/2}.$$

Let $Q(a, t)$ be this region on the unit sphere (i.e., $u = \cos\theta$, $\|v\| = \sin\theta < \frac{a}{t}\left(\frac{1-t^2}{1-a^2}\right)^{1/2}$). From (9.6) we have

(9.8) $$a^{N-1} R^N \frac{\sigma_N}{N} = O(a^N) + \int_a^1 \left[2 \int_{Q(a,t)} \psi(tRx) dA(x)\right] t^{N-1} R^N dt.$$

Divide by $a^{N-1} R^N$ and use (9.5) to find

$$\frac{\sigma_N}{N} = O(a) + \frac{\sigma_{N-1}}{N-1} \int_a^1 2(\psi(tR) + o(1)) \left[(1 - t^2)^{\frac{N-1}{2}} + o(1)\right] dt$$

where $o(1)$'s are as $a \to 0$. We conclude that for all $R > 0$

(9.9) $$\frac{\sigma_N}{2\sigma_{N-1}} \frac{N-1}{N} = \int_0^1 \psi(tR)(1 - t^2)^{\frac{N-1}{2}} dt.$$

Now the left side of (9.9) also equals $\int_0^1 (1 - t^2)^{\frac{N-1}{2}} dt$ as is easily checked. By the Wiener tauberian theorem, applied to the bounded function $\psi(r)$, $r > 0$, and $t(1 - t^2)^{\frac{N-1}{2}} \chi_{(0,1)}(t) \in L^1\left(\mathbf{R}^+, \frac{dt}{t}\right)$, we have

$$\int_0^1 \psi(tR) dt = 1 \quad (R > 0).$$

Since $\psi$ is continuous, $\psi(r) = 1$, $r > 0$. Therefore, $\nu(dx) = \psi(x)m(dx) = m(dx)$, and Theorem 9.1 is proved in this special case.

Now let there be given an arbitrary measure $\nu(dx)$ which satisfies the hypotheses of Theorem 9.1. Let $H = \mathbf{R}^+ \times K$, $K = \mathrm{SO}(N)$. Identify $\mathbf{R}^N \setminus \{0\}$ with $H/\mathrm{SO}(N-1)$ by the map $h \to he_N = tke_N$, $t > 0$, $k \in K$. There is a canonical lift $\mu$ of $\nu$ to a Borel measure on $H$ which is right invariant under (the embedded) $\mathrm{SO}(N-1)$. Because $\nu$ is *even*, $\mu$ is right invariant under $L = \{k \in K \mid k\{\pm e_N\} = \{\pm e_N\}\}$.

Let $\varphi \in C_c(H)$ be such that

$$(9.10) \qquad \int_H \varphi(u,\ell) u^{-N-1} du\, m_K(d\ell) = 1.$$

Use $\varphi$ and $\mu$ to set up a function $\hat\psi$ on $H$, where

$$(9.11) \qquad \hat\psi(s, k_0) = \frac{1}{s^N} \int_H \varphi\left(\frac{s}{t}, kk_0^{-1}\right) \mu(d(t,k)).$$

The right $L$-invariance of $\mu$ implies $\hat\psi$ is right $L$-invariant. Therefore, $\hat\psi$ determines $\psi(\cdot)$ on $\mathbf{R}^N \setminus \{0\}$, where

$$(9.12) \qquad \psi(x) = \hat\psi(s, k_0) \qquad (x = sk_0 e_N).$$

It is evident that $\psi$ is continuous and even.

Let $0 < r_0 < r_1 < \infty$ be such that $\varphi$ is supported on $[r_0, r_1] \times K$. Given $s > 0$, the set of $t$ such that $\varphi\left(\frac{s}{t}, kk_0^{-1}\right) \neq 0$ for some $k \in K$ satisfies $t < s/r_0$, and therefore

$$\begin{aligned} |\hat\psi(s,k_0)| &\leq \|\varphi\|_\infty \frac{\sigma_N}{N} \left(\frac{s}{r_0}\right)^N \cdot \frac{1}{s^N} \\ &= \|\varphi\|_\infty \frac{\sigma_N}{Nr_0^N}. \end{aligned}$$

Fix an ellipsoid $E$ centered at $0$. We have

$$\begin{aligned} &\int_{\mathbf{R}^N} \psi(x)\chi_E(x) m(dx) \\ &= \int_{\mathbf{R}^+ \times K} \psi(sk_0 e_N)\chi_E(ske_N) s^{N-1} ds\, m_K(dk_0) \\ &= \int_{\mathbf{R}^+ \times K} \hat\psi(s,k_0) \chi_E(sk_0 e_N) s^{N-1} ds\, m_K(dk_0) \\ &= \int_{\mathbf{R}^+ \times K} \frac{\chi_E(sk_0 e_N)}{s^N} \int_H \varphi\left(\frac{s}{t}, kk_0^{-1}\right) \mu(d(t,k)) s^{N-1} ds\, m_K(dk_0) \\ &= \int_{\mathbf{R}^+ \times K} \left[\int_H \chi_E(tu\ell^{-1}k) \mu(d(t,k))\right] \varphi(u,\ell) \frac{du}{u}\, m_K(d\ell) \end{aligned}$$

where in the last line we have substituted $u = \frac{s}{t}$ and $\ell = kk_0^{-1}$, $k_0 = \ell^{-1}k$. The expression in brackets has, by assumption, the value $m\left(\frac{\ell E}{u}\right) = \frac{1}{u^N}m(E)$. By (9.10) we have

$$\int_{\mathbf{R}^N} \psi(x)\chi_E(x)m(dx) = m(E).$$

Since $\psi \in C(\mathbf{R}^N \setminus \{0\})$ is bounded, the first part of the argument implies $\psi \equiv 1$. It follows that for any $\varphi \in C_c(H)$

(9.13) $$\int_H \varphi\left(\frac{s}{t}, kk_0^{-1}\right) \mu(d(t,k)) = s^N \int_H \varphi(u, \ell)u^{-N-1}du\, m_K(d\ell).$$

Since (9.13) is also true when $\mu$ is replaced by the measure $t^{N-1}dt m_K(dk)$, and since $\varphi \in C_c(H)$ is arbitrary, it must be that $\mu(d(t,k)) = t^{N-1}dt m_K(dk)$. This implies $\nu = m$ is Lebesgue measure, and Theorem 9.1 is proved.

Theorem 9.1 has been stated for Lebesgue measure since that is the immediate application. A small modification of the two step proof establishes

THEOREM 9.4. *Let $\nu_1, \nu_2 \in \mathcal{M}_N$ be even, and assume $\nu_1(E) = \nu_2(E)$ for every ellipsoid $E$ with center $0$. Then $\nu_1 = \nu_2$.*

*Proof.* Let $\nu = \nu_1 - \nu_2$, and first assume $\nu$ has the form (9.4). One is led by the same argument to the relation (9.9) with $\frac{\sigma_N}{2\sigma_{N-1}}\frac{N-1}{N}$ replaced by 0. One then infers $\nu = 0$. In the general case lift $\nu$ to $\mu$ on $H = R^+ \times K$ as in the paragraph which contains (9.10), and given $\varphi \in C_c(H)$ define $\hat{\psi}$ by (9.11). $\hat{\psi}$ determines $\psi$ by (9.12), and one finds $\psi$ is bounded and even with integral zero over every ellipsoid $E$ centered at 0. Then $\psi = 0$, whence $\hat{\psi} = 0$ and, letting $\varphi$ vary, $\mu = 0$. Details are left to the reader.

## 10. Uniform distribution

Theorem 9.1 and Corollary 9.2 have as an almost immediate consequence a sort of "Weyl criterion" for a notion of uniform distribution on $\mathbf{R}^N$. $B = B(0, 1)$ is the unit ball in $\mathbf{R}^N$.

THEOREM 10.1. *Let $\{\nu_\alpha \mid \alpha \in A\}$ be a net of even, locally finite Borel measures on $\mathbf{R}^N$. Assume there exist $c < \infty$ and a dense set $F \subseteq G \times \mathbf{R}^+$ such that*

(10.2) $$\lim_{\alpha \in A} \nu_\alpha(tgB) = ct^N \frac{\sigma_N}{N} \qquad ((g,t) \in F).$$

*Then*

(10.3) $$\lim_{\alpha \in A} \int_{\mathbf{R}^N} \psi(x)\nu_\alpha(dx) = c \int_{\mathbf{R}^N} \psi(x)m(dx) \qquad (\psi \in C_c(\mathbf{R}^N))$$

*where $m(dx)$ is Lebesgue measure.*

*Proof.* It is easy to see that the hypothesis on $F$ is also true on $\overline{F} = G \times \mathbf{R}^+$. Then (10.3) is an immediate consequence of Corollary 9.2. □

We shall also apply Corollary 9.2 to nets with index set $A = \mathbf{R}^+$ which depend upon a "parameter" $\nu \in \mathcal{M}_N^e$ (Section 8). Let $\mu$ be a Siegel measure supported on $\mathcal{M}_N^e$ and for each $\nu \in \mathcal{M}_N^e$ and $R > 0$ define

$$\nu_R = T_R \nu \tag{10.4}$$

where $T_R \nu(U) = R^{-N} \nu(RU)$, as in Section 8. Theorem 7.1 and the $G$-invariance of $\mu$ imply

$$\lim_{R \to \infty} \left\| g^{-1} \nu_R(B) - c(\mu) \frac{\sigma_N}{N} \right\|_1 = \lim_{R \to \infty} \left\| (g^{-1}\nu)_R(B) - c(\mu) \frac{\sigma_N}{N} \right\|_1 = 0.$$

As $\nu_R(tB) = t^N \nu_{tR}(B)$, we also have

$$\lim_{R \to \infty} \left\| g^{-1} \nu_R(tB) - c(\mu) t^N \frac{\sigma_N}{N} \right\|_1 = 0. \tag{10.5}$$

Let $\psi \in C_c(\mathbf{R}^N)$ and $\epsilon > 0$ be given. Corollary 9.2 implies there exist $\delta > 0$ and $(g_j, t_j) \in G \times \mathbf{R}^+$, $1 \leq j \leq r$, such that whenever $\nu$ is an even Borel measure such that $|g_j^{-1} \nu(t_j B) - c(\mu) t_j^N \frac{\sigma_N}{N}| < \delta$, $1 \leq j \leq r$, then

$$\left| \int \psi(x) \nu(dx) - c(\mu) \int \psi(x) m(dx) \right| < \epsilon. \tag{10.6}$$

From (10.5) we have

$$\lim_{R \to \infty} \sum_{j=1}^r \left\| g_j^{-1} \nu_R(t_j B) - c(\mu) t_j^N \frac{\sigma_N}{N} \right\|_1 = 0.$$

It follows therefore that

$$\lim_{R \to \infty} \hat{\psi}(\nu_R) = c(\mu) \hat{\psi}(m)$$

exists in $\mu$-measure for every $\psi \in C_c(\mathbf{R}^N)$. Since $\hat{\chi}_{tB}(\nu_R)$ converges in $L^1(\mu)$ as $R \to \infty$ for every $t > 0$, it also follows that $\{\hat{\psi}(\nu_R) \mid R > 0\}$ is uniformly integrable for every $\psi \in C_c(\mathbf{R}^N)$. As a consequence we have

THEOREM 10.7. *Let $\mu \in \mathcal{P}(\mathcal{M}_N^e)$ be a Siegel measure, and let $\nu_R$, $R > 0$, $\nu \in \mathcal{M}_N^e$ be as in (10.4). For all $\psi \in C_c(\mathbf{R}^N)$*

$$\lim_{R \to \infty} \left\| \int \psi(x) \nu_R(dx) - c(\mu) \int \psi(x) m(dx) \right\|_1 = 0. \tag{10.8}$$

Let $\mu \in \mathcal{P}(\mathcal{M}_N^e)$ be a Siegel measure, and let $E \subseteq \mathbf{R}^+$ be an unbounded set. Assume it is known to be true that for $\mu$-a.e. $\nu$ there is a dense set $F(\nu) \subseteq G \times \mathbf{R}^+$ such that

$$\lim_{\substack{R \to \infty \\ R \in E}} g^{-1} \nu_R(tB) = c(\mu) t^N \frac{\sigma_N}{N} \qquad ((g,t) \in F(\nu)).$$

Theorem 10.1 implies that if $m(dx)$ is Lebesgue measure, then for $\mu$-a.e. $\nu$

$$\lim_{\substack{R \to \infty \\ R \in E}} \nu_R = c(\mu) m$$

in the $C_c(\mathbf{R}^N)$ topology. Theorems 5.28, 10.7, the proof of Theorem 10.1 and a Borel-Cantelli argument imply

THEOREM 10.9. *Let $\mu \in \mathcal{P}(\mathcal{M}_N^e)$ be a Siegel measure. There exists a sequence $R_n \to \infty$ such that for $\mu$-a.e. $\nu$*

(10.10) $$\lim_{n \to \infty} g^{-1} \nu_{R_n}(B) = c(\mu) \frac{\sigma_N}{N} \qquad (g \in G)$$

(10.11) $$\lim_{n \to \infty} \nu_{R_n} = c(\mu) m.$$

*Convergence in the second line is in the $C_c(\mathbf{R}^N)$ topology.*

## 11. Quadratic differentials

Fix $p$, $n > 0$, and let $M_{p,n}$ be a closed oriented surface ($M_p$) of genus $p$ with $n$ punctures ($S_n$). $H(p,n)$ denotes the group of orientation-preserving homeomorphisms of $M_{p,n}$ with identity component $H_0(p,n)$. Set $\mathrm{Map}(p,n) = H(p,n)/H_0(p,n)$, the mapping class group.

$\Omega^+(p,n)$ denotes the set of admissible positive $F$-structures on $M_{p,n}$. A *positive $F$-structure* is an atlas $u$ on $M_{p,n}$ with three properties: (i) coordinate transitions are locally translations, (ii) $u$ is compatible with orientation and (iii) $u$ is maximal relative to (i) and (ii). The euclidean metric lifts via $u$ charts to a Riemannian flat metric $g(u)$, and $u$ is *admissible* if $M_p$ is the completion of $M_{p,n}$ for the $g(u)$ geodesic function. $u$ determines a complex structure $J(u)$ and nowhere zero holomorphic 1-form $\omega(u)$ ($= f^* dz$ for $u$ chart functions $f$); admissibility is equivalent to the requirement that $J(u)$ extend to $M_p$ and $\omega(u)$ extend as a holomorphic 1-form. If $\omega(u)$ has a zero of order $\nu$ at $s \in S_n$, $g(u)$ has a cone singularity with cone angle $2\pi(\nu + 1)$ at $s$.

Define $\mathcal{M}^+(p,n) = \Omega^+(p,n)/H_0(p,n)$. $\mathcal{M}^+(p,n)$ carries a complete metric with respect to which $\mathrm{Map}(p,n)$ acts properly discontinuously by isometries (cf. [22, Section 1]).

The map $u \to \widehat{\omega(u)} \in H^1_{\mathbf{C}}(M_p, S_n)$ is a local homeomorphism which endows $\mathcal{M}^+(p,n)$ with the structure of a complex manifold ([22, Remark 7.22]). In local coordinates $\mathrm{Map}(p,n)$ is represented by $\mathrm{GL}(2p - 1 + n, \mathbf{Z})$ acting linearly on $H^1_{\mathbf{C}}(M_p, S_n) \cong \mathbf{C}^{2p-1+n}$. Therefore, $\mathrm{Map}(p,n)$ preserves not only the complex structure on $\mathcal{M}^+(p,n)$ but also the lift to $\mathcal{M}^+(p,n)$ of the euclidean volume form on $H^1_{\mathbf{C}}(M_p, S_n)$, made canonical by the requirement that the lattice $H^1_{\mathbf{Z}}(M_p, S_n)$ have covolume one ([22, Theorem 7.17]).

Let $G = \mathrm{SL}(2, \mathbf{R})$. If $g \in G$ and $u \in \Omega^+(p, n)$, then $gu$ is defined by postcomposition of $u$ chart functions with the $\mathbf{R}$-linear transformation $g$. One finds $gu \in \Omega^+(p, n)$ and that $\omega(gu)$ has the same zero structure on $S_n$ as $\omega(u)$. In terms of the common underlying real analytic structure on $M_{p,n}$ we have

$$\omega(gu) = \alpha \omega(u) + \overline{\beta}\,\overline{\omega(u)} \tag{11.1}$$

when $g \in G$ is represented by $\begin{pmatrix} \alpha & \overline{\beta} \\ \beta & \overline{\alpha} \end{pmatrix} \in \mathrm{SU}(1,1)$.

We find therefore that

$$\widehat{\omega(gu)} = \alpha \widehat{\omega(u)} + \overline{\beta}\,\overline{\widehat{\omega(u)}}. \tag{11.2}$$

The $G$ action on $\Omega^+(p, n)$ commutes with the action of $H_0(p, n)$ and therefore descends to $\mathcal{M}^+(p, n)$ where it is, by (11.2), real analytic. Since $\det g = 1$ implies (11.2) is euclidean volume preserving $G$ preserves volume on $\mathcal{M}^+(p, n)$.

If $m = [u] \in \mathcal{M}^+(p, n)$, define $V(m) = \frac{i}{2} \int_{M_p} \omega(u) \wedge \overline{\omega(u)}$. $V(\cdot)$ is real analytic without critical points, and therefore $\mathcal{M}_1^+(p, n) = V^{-1} 1$ is a real analytic real codimension one submanifold. $V(\cdot)$ is $G$-invariant, and therefore $dV$ and the canonical volume element on $\mathcal{M}^+(p, n)$ determine a canonical $G$-invariant volume element on $\mathcal{M}_1^+(p, n)$. This volume element satisfies

$$\mathrm{Vol}\left(\mathcal{M}_1^+(p, n) / \mathrm{Map}(p, n)\right) < \infty. \tag{11.3}$$

(See [12], [21], [22], [15].)

In what follows $\mathcal{M}$ denotes a fixed connected component of $\mathcal{M}_1^+(p, n)/\mathrm{Map}(p, n)$. $\lambda$ denotes the $G$-invariant probability measure obtained, using (11.3), by normalizing the natural image measure on $\mathcal{M}$. We recall that $(\mathcal{M}, G, \lambda)$ is *ergodic* ([12], [21], [26]).

*Remark* 11.4. Let $(\mathcal{M}, G)$ be as above. We recall that if $m \in \mathcal{M}$ and $\Gamma(m) = \{g \in G \mid gm = m\}$, then $\Gamma(m)$ is a discrete subgroup. For a dense set of $m$ $\Gamma(m)$ is a lattice. (For example, $\Gamma(m)$ is a lattice if $\widehat{\omega(u)}$ ($m = [u]\,\mathrm{Map}(p, n)$) is projectively a rational class. [23]) For such $m$ the normalized Haar measure on $G/\Gamma(m)$ determines an ergodic invariant probability measure on the orbit $Gm \subseteq \mathcal{M}$. Another measure of interest arises when there exists $\tau \in H(p, n)$ such that $\tau^2 = \mathrm{Id}$ and $\mathrm{Fix}\,\tau \subseteq S_n$. Define $\mathcal{M}_\tau^+(p, n)$ to be the set of $[u]$ such that for some $\varphi \in H_0(p, n)$ and $\tau_\varphi = \varphi \tau \varphi^{-1}$, $\tau_\varphi^* \omega(u) = -\omega(u)$. $\mathcal{M}_\tau^+(p, n)$ is a closed set and complex submanifold which also carries a natural $G$-invariant volume. Proceeding by analogy with the discussion above, one finds that $\mathcal{M}_1^+(p, n) \cap \mathcal{M}_\tau^+(p, n)$ also carries a $G$-invariant volume and the projection in $\mathcal{M}_1^+(p, n)/\mathrm{Map}(p, n)$ has finite total volume. (The projection depend only upon $[\tau] \in \mathrm{Map}(p, n)$, of course.)

If $\mathcal{M}_{[\tau]}$ is a component of $\left(\mathcal{M}_1^+(p, n) \cap \mathcal{M}_\tau^+(p, n)\right)/\mathrm{Map}(p, n)$ then $\mathcal{M}_{[\tau]}$ carries a natural normalized $G$-invariant ergodic measure ([26, Theorem 6.14]).

Of course, $\mathcal{M}_{[\tau]} \subseteq \mathcal{M}$ for some component $\mathcal{M} \subseteq \mathcal{M}_1^+(p,n)/\operatorname{Map}(p,n)$. Examples such as $\mathcal{M}_{[\tau]}$ arise from lifting via 2-sheeted branched coverings quadratic differentials in genus $p' \geq 0$ to holomorphic 1-forms which are odd relative to the (branched-) covering transformation.

## 12. Siegel measures and quadratic differentials

Let $\mathcal{S}(p,n)$ be the set of free homotopy classes of simple closed curves in $M_{p,n}$. If $m = [u] \in \mathcal{M}^+(p,n)$, we define $\mathcal{S}(m)$ to be the set of $\gamma \in \mathcal{S}(p,n)$ such that $\gamma$ has a closed $g(u)$-geodesic representative for any $u \in m$. When $\gamma \in \mathcal{S}(m)$, there is a number $a(m,\gamma) > 0$ which is for any $u \in m$ the area of the cylinder of closed $g(u)$ geodesics which represent $\gamma$. There is also a symmetric pair of vectors $\pm v(m,\gamma)$ giving the length and possible directions, determined by any atlas $u \in m$, of closed geodesics which represent $\gamma$. It is an elementary consequence of the definition of the metric on $\mathcal{M}^+(p,n)$ ([22, Section 1]) that for any $\gamma \in \mathcal{S}(p,n)$ and $s \geq 0$ the set

$$(12.1) \qquad U(\gamma, s) = \{ m \in \mathcal{M}^+(p,n) \mid \gamma \in \mathcal{S}(m), a(m,\gamma) > s \}$$

is open and the pair $\pm v(m,\gamma)$ varies continuously on $U(\gamma, s)$. In particular, if $\psi \geq 0$ is a Borel function on $\mathbf{R}^2$, the function

$$(12.2) \qquad T_s \psi(m) = \sum_{\gamma \in \mathcal{S}(p,n)} \chi_{U(\gamma,s)}(m) \psi(\pm v(m,\gamma))$$

is Borel. (The $\pm$ indicates two summands for each $\gamma \in \mathcal{S}(p,n)$.)

A starting point for the present work has been the theorem of Masur which is cited in the introduction. If we define $N(m, s, R)$ to be the growth function of

$$(12.3) \qquad \Pi(m,s) = \{ \pm v(m,\gamma) \mid \gamma \in \mathcal{S}(m), m \in U(\gamma, s) \}$$

then

$$(12.4) \qquad N(m, 0, R) = O(R^2) \qquad (R \to \infty).$$

Masur also establishes a lower quadratic bound which is not necessary for the present discussion. If $V(m) = 1$, and if $0 < s < V(m)$, then as is implicit in [13] there is a uniform constant $C(s, p, n) < \infty$ such that

$$(12.5) \qquad N(m, s, R) < C(s, p, n)(R^2 + 1) \qquad (V(m) = 1, \ R > 0).$$

This implies that if $\psi \geq 0$ is bounded, Borel with compact support, then $T_s \psi$ ((12.2)) is uniformly bounded on $\mathcal{M}_1^+(p,n)$.

Now let us suppose given an ergodic $G$-invariant Borel probability measure $\eta$ on $\mathcal{M}_1^+(p,n)/\operatorname{Map}(p,n)$. Let $0 \leq s < 1$, and let $\nu_{m,s}$ be counting measure on the set $\Pi(m, s)$ in (12.3). From the definition of $G$ action in Section 11 it follows for any $m$, $s$ that

(12.6) $$g\Pi(m,s) = \Pi(gm,s) \qquad (g \in G).$$

If we couple this with (12.4)–(12.5), we find

(12.7) $$\nu_{m,s} \in \mathcal{M}_2$$
$$\nu_{gm,s} = g\nu_{m,s}.$$

The discussion of (12.2) implies the map

(12.8) $$\xi_s(m) = \nu_{m,s}$$

is Borel; therefore the measures

(12.9) $$\eta_s = \xi_s \eta \qquad (0 \leq s < 1)$$

are Borel measures on $\mathcal{M}_2$. Since $\eta$ is $G$-invariant and ergodic, (12.7) implies $\eta_s$ is $G$-invariant and ergodic. We have

THEOREM 12.10. *Let $\eta$ be a $G$-invariant ergodic probability measure on $\mathcal{M}_1^+(p,n)/\operatorname{Map}(p,n)$. For all $s$ such that $0 \leq s < 1$ the measure $\eta_s = \xi_s \eta$, defined by (12.8)–(12.9) is a Siegel measure.*

Our main result concerning quadratic differentials, Theorem 12.11 below, is now a corollary of Theorems 12.10, 6.5, 7.1 and 10.6:

THEOREM 12.11. *Let $\eta$ be an ergodic $G$-invariant Borel probability measure on $\mathcal{M}_1^+(p,n)/\operatorname{Map}(p,n)$. There exist constants $c(\eta,s) < \infty$, $0 \leq s < 1$, such that the following statements obtain*:

I. *If $\psi \geq 0$ is Borel on $\mathbf{R}^2$, and if $T_s \psi$, $0 \leq s < 1$, is defined by (12.2), then*

(12.12) $$\int_{\mathcal{M}_1^+(p,n)/\operatorname{Map}(p,n)} T_s \psi(m) \eta(dm) = c(\eta,s) \int_{\mathbf{R}^2} \psi(x) dx.$$

II. *If $0 \leq s < 1$, and if $N(m,s,R)$ is defined as the growth function of $\Pi(m,s)$ in (12.3), then*

(12.13) $$\lim_{R \to \infty} \frac{N(m,s,R)}{R^2} = c(\eta,s)\pi$$

*in $L^1\left(\mathcal{M}_1^+(p,n)/\operatorname{Map}(p,n), \eta\right)$.*

III. *If $0 \leq s < 1$, if $\psi \in C_c(\mathbf{R}^2)$ and if $\psi_R(v) = \psi\left(\frac{v}{R}\right)$, then*

(12.14) $$\lim_{R\to\infty} \frac{1}{R^2} (T_s \psi_R)(m) = c(\eta, s) \int_{\mathbf{R}^2} \psi(x) dx$$

in $L^1\left(\mathcal{M}_1^+(p,n)/\operatorname{Map}(p,n), \eta\right)$.

*Remark* 12.15. While the focus in this paper has been on periodic trajectories, entirely analogous results follow by the same techniques for sets in the plane which represent simple geodesics joining cone points (i.e. zeros) for the metrics $|\omega|^2$, $\omega$ a holomorphic 1-form. The requisite quadratic (upper) estimate is also due to Masur ([13]).

## 13. Properties of $c(\eta, s)$

Notations are as in Section 12. If $\psi \geq 0$ on $\mathbf{R}^2$, the monotone convergence theorem, applied in (12.2) to counting measure on $\mathcal{S}(p,n)$, implies the function $s \to T_s\psi(m)$ is for each $m$ continuous from the right on $[0,1)$. If $\psi$ is assumed to be Borel with finite positive integral over $\mathbf{R}^2$, (12.12) and the monotone convergence theorem imply $c(\eta, \cdot)$ is also continuous from the right on $[0,1)$.

It is not true in general that $c(\eta, \cdot) \in C([0,1))$. When $\eta$ is concentrated on an orbit (see the first part of Remark 11.4), the range of $a(m, \gamma)$ is a finite set ([23]) and $c(\eta, \cdot)$ is a step function which is, in general, not constant. However, the fact that $T_0\psi \in L^1(\eta)$, $\psi$ as above, and the dominated convergence theorem imply that if $s \to T_s\psi(m)$ is for each $s_0 \in (0,1)$ and a.e. $[m]$ left continuous at $s_0$, then $c(\eta, \cdot) \in C([0,1))$:

PROPOSITION 13.1. *Let $\eta \in \mathcal{P}\left(\mathcal{M}_1^+(p,n)/\operatorname{Map}(p,n)\right)$ be invariant and ergodic. If $s_0 \in (0,1)$ is such that*

$$\eta\left\{[m] \mid a(m, \gamma) = s_0 \text{ for some } \gamma \in \mathcal{S}(m)\right\} = 0$$

*then $c(\eta, \cdot)$ is continuous at $s_0$.*

In what follows $\mathcal{M}$ denotes a fixed topological component of $\mathcal{M}_1^+(p,n)/\operatorname{Map}(p,n)$ and $\lambda$ the $G$-invariant probability measure obtained, using (11.3), from normalizing the natural image measure on $\mathcal{M}$. We recall that $(\mathcal{M}, G, \lambda)$ is *ergodic* ([12], [21], [26]).

The discussion which follows is local. Therefore we fix $\mathcal{M}$ and $[m_0] \in \mathcal{M}$ and work with $m_0$ and a fixed $\gamma \in \mathcal{S}(m_0)$ such that $a(m_0, \gamma) = s_0$. By definition $V(m_0) = 1$. Let $U(m_0)$ be an open set in $\mathcal{M}_1^+(p,n)$ containing $m_0$ and with the properties (i) $\gamma \in \mathcal{S}(m)$, $m \in U(m_0)$ and (ii) $\hat{\omega}(\cdot)$ is schlicht on $U(m_0)$. The functions $V(\cdot)$ and $a(\cdot, \gamma)$ are quadratic forms in the coordinate (ii). In particular, if $V(m) = 1$ implies $a(m, \gamma) = s_0$ in this coordinate, then $a(\cdot, \gamma) - s_0 V(\cdot)$ is identically zero on $U(m_0)$. We shall observe this implies $s_0 = 1$ and $2p - 1 + n = 2$, i.e., $p = 1 = n$.

With notations as above we consider separately the cases $s_0 = 1$ and $s_0 < 1$.

*Case* 1. $s_0 = 1$. Let $C(m_0, \gamma)$ be the cylinder which corresponds to $\gamma$. Since $a(m_0, \gamma) = 1 = V(m_0)$, there is a parallelogram $P$ of area one and a gluing-by-translation rule on $\partial P$ such that $P/\sim$ equipped with its natural 1-form ('$dz$') realizes $m_0$. The gluing rule is pure translation between one pair of parallel edges and piecewise translation between another pair of edges. If the latter gluing is not pure translation, that is, if $2p - 1 + n > 2$, it is clear by inspection that $a(m, \gamma) \not\equiv 1$ on $U(m_0) \cap \{\|m\| = 1\}$. If $p = 1 = n$, then $\mathcal{M}$ is an orbit.

*Case* 2. $0 < s_0 < 1$. Geodesic triangulations of $C(m_0, \gamma)$ and $M_{p,n} \setminus C(m_0, \gamma)$ can be used to define two nonempty sets of geodesics, $A$ and $B$, joining points of $S_n$, such that $A \cap B = \emptyset$, $A \cup B$ span $H_1(M_p, S_n)$ and the areas of $M_{p,n} \setminus C(m_0, \gamma)$ and $C(m_0, \gamma)$ are quadratic forms $Q_0(\hat{\omega}(\cdot))$ and $Q_1(\hat{\omega}(\cdot))$, with $Q_0$ depending upon $\hat{\omega}(m)\big|_A$ and $Q_1$ depending upon $\hat{\omega}(m)\big|_{A \cup B}$. It is possible to vary $m$ in $U(m_0)$ in such a way that $\hat{\omega}(\cdot)\big|_A$ remains constant while $Q_1(\hat{\omega}(\cdot))$ does not. We have by assumption $a(\cdot, \gamma) = s_0 V(\cdot)$ on $U(m_0)$, and therefore

$$(13.2) \qquad Q_1(\hat{\omega}(m)) = s_0 \left( Q_0(\hat{\omega}(m)) + Q_1(\hat{\omega}(m)) \right).$$

Varying $m$ as above we find that $s_0 = 1$, $Q_0 \equiv 0$, a contradiction.

THEOREM 13.3. *Let $2p - 1 + n > 2$, and let $\mathcal{M}$ be a component of $\mathcal{M}_1^+(p, n)/\operatorname{Map}(p, n)$ equipped with its invariant normalized volume $\lambda$. Then $c(\lambda, \cdot) \in C([0, 1))$. Moreover, $c(\lambda, \cdot)$ is strictly decreasing on $[0, 1)$ and $c(\lambda, 1^-) = 0$.*

*Proof.* Case 1 above implies that when $2p - 1 + n > 2$, then for $\lambda$-a.e. $[m] \in \mathcal{M}$ the function $s \to T_s \psi(m)$ vanishes as $s \to 1$ as soon as there exists $s$ such that $T_s \psi(m) < \infty$. Assuming $\psi \geq 0$ is integrable over $\mathbf{R}^2$, this latter requirement is satisfied for $\lambda$-a.e. $[m]$. Therefore $c(\lambda, 1^-) = 0$. To prove that $c(\lambda, \cdot)$ is strictly decreasing it is sufficient to prove there exists $[m_s] \in \mathcal{M}$, $0 < s < 1$, and $\gamma \in \mathcal{S}(m_s)$ such that $a(m_s, \gamma) = s$. To this end fix any $m \in \mathcal{M}_1^+(p, n)$ such that $[m]$ projects to a point of the given component $\mathcal{M}$. Choose any cylinder of closed geodesics for $m$, and observe that this cylinder may be elongated or shortened so as to occupy as large or small a relative portion of the total volume of (the altered) $m$. Normalizing the altered $m$ produces $m_s$ for any $s \in (0, 1)$. □

*Question* 13.4. If $2p - 1 + n > 2$, and if $(\mathcal{M}, \lambda)$ are as in the theorem, is there a simple formula for $b(\lambda, s) = \frac{c(\lambda, s)}{c(\lambda, 0)}$?

## 14. Pointwise statements

Let $G = \mathrm{SL}(N, \mathbf{R})$, $K = \mathrm{SO}(N)$, and let $A_N^+$ be as in Section 5. $\pi$ denotes a continuous unitary representation of $G$ on a Hilbert space $H$ such that (a) $\pi$ admits a cyclic vector $u$, $\|u\| = 1$, which is fixed by $K$ and (b) $\pi$ admits no nonzero invariant vectors. In addition we assume (c) if $N = 2$, then $\pi$ does not almost have invariant vectors; that is, there exist $\epsilon > 0$ and a compact set $C \subseteq G$ such that

$$(14.1) \qquad \underset{c \in C}{\mathrm{Max}}\, \|\pi(c)v - v\| \geq \epsilon \|v\| \qquad (v \in H).$$

If $g \in G$, express $g$ as $g = k_1(g) a^+(g) k_2(g)$ with $k_j(g) \in K$, $j = 1, 2$, and $a^+(g) \in A_N^+$. Define $\sigma(g)$ to be the minimum of the ratios between diagonal entries of $a^+(g) = \mathrm{diag}(a_1(g), \ldots, a_N(g))$, i.e., $\sigma(g) = a_N(g)/a_1(g)$. The assumption (a)–(c) above together with estimates in [9, Chapter V], imply there exists $\eta > 0$ such that if $\sigma(g)$ is sufficiently small, then

$$(14.2) \qquad |\langle \pi(g)u, u \rangle| \leq \sigma(g)^\eta \qquad (g \in G)$$

where $\langle \cdot, \cdot \rangle$ is the inner product on $H$.

We specialize $g$ in what follows. We define

$$a(t) = \mathrm{diag}\left(e^{(N-1)t}, e^{(N-3)t}, \ldots, e^{(1-N)t}\right) \in A_N^+, \quad t > 0.$$

We have

$$(14.3) \qquad \|Q_\pi \circ \pi(a(t))u\|^2 = \int_K \langle \pi(a(t)^{-1} k a(t))u, u \rangle\, m_K(dk).$$

To estimate the size of the integrand in (14.3) in terms of (14.2) it is necessary to estimate the first and last diagonal entries of $a^+(a^{-1}(t) k a(t))$, $k \in K$. The first diagonal entry, denoted $a_1^+$, satisfies

$$N^{1/2} a_1^+ \geq \|a^{-1}(t) k a(t)\|_{\mathrm{HS}}$$

where $\|\cdot\|_{\mathrm{HS}}$ is Hilbert-Schmidt norm. If $S(k) = \sum_{j=1}^N (1 - k_{jj}^2)$, the Hilbert-Schmidt norm satisfies

$$\|a^{-1}(t) k a(t)\|_{\mathrm{HS}} \geq \mathrm{Max}\left((N - S(k))^{1/2}, \beta e^{2t} S(k)^{1/2}\right)$$

where $\beta > 0$ is a dimensional constant. Since $S(k) = S(k^{-1})$, the last diagonal entry, $a_N^+$, of $a^{-1}(t) k a(t)$ satisfies

$$N^{1/2}(1/a_N^+) \geq \|a^{-1}(t) k^{-1} a(t)\|_{\mathrm{HS}} \geq \mathrm{Max}\left((N - S(k))^{1/2}, \beta e^{2t} S(k)^{1/2}\right)$$

and therefore

$$\sigma\left(a^{-1}(t) k a(t)\right) \leq \mathrm{Min}\left(\frac{N}{N - S(k)}, \frac{N}{\beta^2 S(k) e^{4t}}\right).$$

Let $\gamma(t) > 0$, to be determined later. The open set $\{k \mid S(k) < \gamma(t)\}$ has Haar measure commensurable with $\gamma(t)^{\binom{N}{2}}$ since $\dim K = \binom{N}{2}$. Let $\eta$ be as in (14.2), and choose $\gamma(t)$ to satisfy

$$\gamma(t)^{\binom{N}{2}} = \left(\frac{1}{\beta^2 \gamma(t) e^{4t}}\right)^{\eta}$$

or

$$\gamma(t) = \left(\beta^2 e^{4t}\right)^{-\frac{\eta}{\eta + \binom{N}{2}}}.$$

Now divide the integral (14.3) into two integrals according to whether $S(k) < \gamma(t)$ or $S(k) > \gamma(t)$. The integral is, by the choice of $\gamma(t)$, bounded by

$$(14.4) \qquad \|Q_\pi \circ \pi(a(t))u\|^2 \leq C e^{-2\xi t}$$

$$\xi = \frac{2\eta \binom{N}{2}}{\eta + \binom{N}{2}}.$$

We now suppose $\mu$ is a Siegel measure on $\mathcal{M}_N$, $N > 1$. It is necessary to assume

$$(14.5) \qquad \widehat{\chi_B} \in L^2(\mu) \qquad (B = B(0, 1)).$$

Define $u(\cdot) \in L^2(\mu)$ by

$$(14.6) \qquad u(\nu) = \widehat{\chi_B}(\nu) - c(\mu)\frac{\sigma_N}{N}.$$

Theorem 6.5 implies $u$ has integral zero. Since $\mu$ is by assumption ergodic, the cyclic subspace $H(u) \subseteq L^2(\mu)$ generated by the $G$-orbit of $u$ contains no invariant vector. If $N = 2$ we assume

*Assumption* 14.7. If $N = 2$, then $H(u)$ does not almost have invariant vectors.

Of course, Assumption 14.7 is the same as (c) in the first paragraph of this section applied to the Siegel measure setting.

With notations as above we apply (14.4) to obtain

$$\int_{\mathcal{M}_N} \left(\int_K \left(\widehat{\chi_B}(a(t)k\nu) - c(\mu)\frac{\sigma_N}{N}\right) m_K(dk)\right)^2 \mu(d\nu) = O(e^{-2\xi t}).$$

Fix $\delta > 1/\xi$, and define $t_n = \delta \log(n+1)$, $n > 0$. Since $e^{-\xi t_n} = 1/(n+1)^{\xi\delta}$ is summable, the Borel-Cantelli lemma implies

$$(14.8) \qquad \lim_{n \to \infty} \int_K \widehat{\chi}_B(a(t_n)k\nu) m_K(dk) = c(\mu) \frac{\sigma_N}{N} \qquad (\mu - \text{a.e. } \nu).$$

Since $a_N(t_n)/a_{N-1}(t_n) = e^{-2t_n} \to 0$, (14.8) and (5.13) imply
(14.9)
$$\lim_{n \to \infty} \int_0^1 \frac{N_\nu\left(\frac{\tau}{a_N(t_n)}\right)}{\left(\frac{\tau}{a_N(t_n)}\right)^N} \left(\frac{2\sigma_{N-1}}{\sigma_N}(1-\tau^2)^{\frac{N-3}{2}}\right) d\tau = c(\mu) \frac{\sigma_N}{N} \qquad (\mu - \text{a.e. } \nu).$$

LEMMA 14.10. *Let $\lambda$ be a Borel measure on $\mathbf{R}^+$, and let $\varphi > 0$ on $\mathbf{R}^+$ be such that for some $\alpha \in \mathbf{R}$ the function $\varphi(t)t^\alpha$ is monotone nondecreasing. If $T_n \nearrow \infty$ in such a way that $T_n/T_{n+1} \to 1$, and if $\lim_{n \to \infty} \psi(T_n) = \ell$ exists, where $\psi(T) = \int_{\mathbf{R}^+} \varphi(T\tau) \lambda(d\tau)$, then $\lim_{T \to \infty} \psi(T) = \ell$.*

*Proof.* For each $T \gg 0$ define $n$ by $T_n \leq T < T_{n+1}$. The assumption on $\varphi$ implies $\psi(t)t^\alpha$ is monotone nondecreasing, and therefore

$$\left(\frac{T_n}{T}\right)^\alpha \psi(T_n) \leq \psi(T) \leq \left(\frac{T_{n+1}}{T}\right)^\alpha \psi(T_{n+1}).$$

Since $T_n/T_{n+1} \to 1$ by assumption, $\lim_{T \to \infty} \psi(T) = \lim_{n \to \infty} \psi(T_n)$ as claimed. □

Collecting results the Wiener Tauberian theorem and Lemma 5.19 imply

THEOREM 14.11. *Let $\mu$ be a Siegel measure, and assume of $\mu$ that $\widehat{\chi}_B \in L^2(\mu)$. If $N > 2$, or if $N = 2$ and Assumption 14.7 is true, then for $\mu$-almost all $\nu$*

$$(14.12) \qquad \lim_{R \to \infty} \frac{N_\nu(R)}{R^N} = c(\mu) \frac{\sigma_N}{N}.$$

*Moreover, for $\mu$-almost all $\nu$ if $\psi \in C_c(\mathbf{R}^N)$*

$$(14.13) \qquad \lim_{R \to \infty} \frac{1}{R^N} \int_{\mathbf{R}^N} \psi\left(\frac{x}{R}\right) \nu(dx) = c(\mu) \int_{\mathbf{R}^N} \psi(y) dy.$$

*Remark* 14.14. Let $\lambda$ be the normalized $G$-invariant volume element on a component $\mathcal{M}$ of $\mathcal{M}_1^+(p,n)/\text{Map}(p,n)$. Let $H_0$ be the orthocomplement of the constants in $L^2(\lambda)$. Should it be the case that the representation $(G, H_0)$ does not almost have invariant vectors, then Theorem 14.11 applies to Parts II and III of Theorem 12.11, at least for $0 < s < 1$. The reason is that if $\mu_s$ is the Siegel measure on $\mathcal{M}_2$ determined by $\lambda$, then $\widehat{\chi}_B \in L^\infty(\mu_s) \subseteq L^2(\mu_s)$, $0 < s < 1$. In view of the fact that $c(\lambda, \cdot) \in C([0,1))$ by Theorem 13.3, it is possible that a pointwise a.e. result valid for $s > 0$ would imply a similar result for $s = 0$. In this regard we raise the

*Question* 14.15. Let $m \in \mathcal{M}_1^+(p,n)$ and $\epsilon > 0$. Does there exist $s = s(m,\epsilon) > 0$ such that

$$N(m,0,R) - N(m,s,R) < \epsilon R^2 \tag{14.16}$$

for large $R$?

## 15. Regular points

Let $\mathcal{M}$ be a component of $\mathcal{M}_1^+(p,n)/\operatorname{Map}(p,n)$. If $\xi \in \mathcal{M}$, define $\mu_\xi \in \mathcal{P}(\mathcal{M})$ by

$$\int_\mathcal{M} \psi(y)\mu_\xi(dy) = \int_K \psi(k\xi) m_K(dk). \tag{15.1}$$

The analysis in [10] may be seen to imply the orbit $G\mu_\xi$ is relatively compact in $\mathcal{P}(\mathcal{M})$ with the $C_c(\mathcal{M})$ topology.

*Definition* 15.2. $\xi \in \mathcal{M}$ shall be called a *regular point* if

$$\lim_{g \to \infty} g\mu_\xi = \eta_\xi \tag{15.3}$$

exists in the $C_c(\mathcal{M})$ topology.

*Example* 15.4. Let $\xi \in \mathcal{M}$ be such that the isotropy group $\Gamma(\xi) = \{g \in G \mid g\xi = \xi\}$ is a lattice in $G$. If $m_{G/\Gamma}$ is normalized Haar measure on $G/\Gamma$, and if $\eta_{G/\Gamma} \in \mathcal{P}(G\xi) \subseteq \mathcal{P}(\mathcal{M})$ is the image of $m_{G/\Gamma}$ under the map $g\Gamma \to g\xi$, then by Theorem 1.2 of [3] $\xi$ is regular and $\eta_\xi = \eta_{G/\Gamma}$.

If $\xi \in \mathcal{M}$ is regular, then $\eta_\xi$ is $G$-invariant and, as noted above, $\eta_\xi \in \mathcal{P}(\mathcal{M})$. $\eta_\xi$ is not *a priori* ergodic, but consideration of (a) the ergodic decomposition of $\eta_\xi$, (b) the fact $T_s\psi$ is bounded for each $s \in (0,1)$ and $\psi \in C_c(\mathbf{R}^2)$ and (c) Theorem 6.5 implies

PROPOSITION 15.5. *Let $\xi \in \mathcal{M}$ be regular. For every $s \in (0,1)$ there exists $c(\xi,s) < \infty$ such that*

$$\int_\mathcal{M} T_s\psi(y)\eta_\xi(dy) = c(\xi,s)\int_{\mathbf{R}^2}\psi(u)du \qquad (\psi \in C_c(\mathbf{R}^2),\ 0 < s < 1). \tag{15.6}$$

If $\xi \in \mathcal{M}$ is a regular point, then because mass is preserved in the limit (15.3), this limit exists in a stronger sense. More precisely, let $C_b(\mathcal{M},\xi)$ be the space of bounded Borel functions on $\mathcal{M}$ which are continuous $\eta_\xi$-a.e. We have

LEMMA 15.7. *If $\xi \in \mathcal{M}$ is regular, the limit (15.3) exists also in the $C_b(\mathcal{M},\xi)$ topology.*

*Proof.* Fix $f \in C_b(\mathcal{M}, \xi)$ and $\epsilon > 0$. Let $Q$ be a compact set in $\mathcal{M}$ such that (a) $\eta_\xi(\mathcal{M}\backslash Q) < \epsilon$ and (b) each $q \in Q$ is a point of continuity of $f$. Let $F$ be a Tietze extension of $f\big|_Q$ such that $F \in C_c(\mathcal{M})$ ($\|F\|_\infty \leq \|f\|_\infty$). For each $q \in Q$ let $U(q, \epsilon)$ be a relatively compact open neighborhood of $q$ such that $|F(q') - F(q)| + |f(q') - f(q)| < \epsilon$, $q' \in U(q, \epsilon)$. Choose a finite set $q_1, \ldots, q_n$ such that $Q \subseteq U(\epsilon) = \bigcup_{j=1}^n U(q_j, \epsilon)$. By construction $|F(q') - f(q')| < 2\epsilon$, $q' \in U(\epsilon)$. Since mass is preserved in (15.3), there exists a compact set $L \subseteq G$ such that $g\mu_\xi(\mathcal{M}\backslash U(\epsilon)) < 2\epsilon$, $g \notin L$. Use $\langle \cdot, \cdot \rangle$ to denote pairing of functions and measures. We have for $g \notin L$

$$\begin{aligned}
&\big|\langle f, g\mu_\xi\rangle - \langle f, \eta_\xi\rangle\big| \\
&\leq \big|\langle f - F, g\mu_\xi\rangle\big| + \big|\langle F, g\mu_\xi\rangle - \langle F, \eta_\xi\rangle\big| + \big|\langle F, \eta_\xi\rangle - \langle f, \eta_\xi\rangle\big| \\
&< (1 + 2\cdot 2\|f\|_\infty)\epsilon + \big|\langle F, g\mu_\xi\rangle - \langle F, \eta_\xi\rangle\big| + (1 + 2\|f\|_\infty)\epsilon.
\end{aligned}$$

Since the second summand on the right converges to 0 as $g \to \infty$, and since $\epsilon > 0$ is arbitrary, it follows that $\lim_{g\to\infty} \langle f, g\mu_\xi\rangle = \langle f, \eta_\xi\rangle$, $f \in C_b(\mathcal{M}, \xi)$, as claimed. $\square$

If $\psi \in C_c(\mathbf{R}^2)$, and if $0 < s < 1$, then $T_s\psi$ is a bounded Borel function on $\mathcal{M}$. $T_s\psi$ is continuous at any $y$ which has no maximal cylinder of closed geodesics of area $s$. It follows that if $\xi$ is regular, and if $s$ is not a point of discontinuity of $c(\xi, \cdot)$, then $T_s\psi \in C_b(\mathcal{M}, \xi)$.

In what follows if $\xi \in \mathcal{M}$ is a regular point, $\Delta(\xi)$ will denote the set of discontinuities of $c(\xi, \cdot)$ in $(0, 1)$. We have

LEMMA 15.8. *If $\xi \in \mathcal{M}$ is a regular point, then*

$$\begin{aligned}
(15.9) \quad \lim_{g\to\infty} \int_{\mathcal{M}} T_s\psi(gy)\mu_\xi(dy) \\
= c(\xi, s) \int_{\mathbf{R}^2} \psi(u)du \qquad \left(\psi \in C_c(\mathbf{R}^2),\ s \in (0,1)\backslash\Delta(\xi)\right).
\end{aligned}$$

Let $B = B(0, 1) \subseteq \mathbf{R}^2$. It is an elementary consequence of (15.9) that the same relation (15.9) also holds for the function $\psi = \chi_B$. Apply Theorems 5.19 and 10.1 to conclude

THEOREM 15.10. *Let $\xi \in \mathcal{M}$ be a regular point. With all notations as above we have*

$$(15.11) \quad \lim_{R\to\infty} \frac{N(g^{-1}\xi, s, R)}{R^2} = c(\xi, s)\pi \qquad (g \in G,\ s \in (0,1)\backslash\Delta(\xi)).$$

*Moreover, if $\psi \in C_c(\mathbf{R}^2)$, and if we set $\psi_R(u) = \psi\left(\frac{u}{R}\right)$, then*

$$(15.12) \quad \lim_{R\to\infty} \frac{1}{R^2} T_s\psi_R(\xi) = c(\xi, s) \int_{\mathbf{R}^2} \psi(u)du \qquad (s \in (0,1)\backslash\Delta(\xi)).$$

Theorems 0.3 and 0.18 of the introduction are both consequences of Theorem 15.10 as it applies in Example 15.4. It is only necessary to recall from [23] that if $\xi \in \mathcal{M}$ is such that $\Gamma(\xi)$ is a lattice, then there is a *finite* set $E(\xi) \subseteq (0,1)$ such that every maximal cylinder for $\xi$ has area $s$ for some $s \in E(\xi)$. This implies $c(\xi, s) = c(\xi, 0^+)$, for $s$ small, and therefore (15.11)–(15.12) obtain also for $s = 0$ and $c(\xi, 0) \stackrel{\text{def}}{=} c(\xi, 0^+)$.

## 16. Nonuniform lattices in $G = \mathrm{SL}(2, \mathbf{R})$

The approach of Section 15 will be used in this section for two purposes. The first is to establish an analog of the combined Theorems 0.3 and 0.18 for an arbitrary nonuniform lattice:

THEOREM 16.1. *Let $\Gamma$ be a nonuniform lattice in $G = \mathrm{SL}(2, \mathbf{R})$, and assume $-I \in \Gamma$. Let $\Lambda$ be a maximal unipotent subgroup of $\Gamma$, and let $v \in \mathbf{R}^2 \setminus \{0\}$ be such that $\Lambda v = v$. There exists a positive, finite constant $c(\Gamma, v)$ such that*

$$(16.2) \qquad \lim_{R \to \infty} \frac{\mathrm{Card}(g\Gamma v \cap B(0, R))}{R^2} = c(\Gamma, v)\pi \qquad (g \in G).$$

*Moreover,*

$$(16.3) \qquad \lim_{R \to \infty} \frac{1}{R^2} \sum_{w \in \Gamma v} \psi\left(\frac{w}{R}\right) = c(\Gamma, v) \int_{\mathbf{R}^2} \psi(u) du \qquad (\psi \in C_c(\mathbf{R}^2)).$$

The proof of Theorem 16.1 will be modelled on the proof of Theorem 15.10. Given the Weyl criterion, Theorem 10.1, the critical issue in Theorem 16.1 is the relation (16.2). As with (0.4) one may prove (16.2) using the theory of Eisenstein series (for $(\Gamma, \Lambda)$) and the Ikehara tauberian theorem, as in [23]. The latter approach also yields an explicit expression for $c(\Gamma, v)$. Therefore, a second purpose of this section will be to observe that $c(\Gamma, v)$ may be computed without the theory of Eisenstein series. In particular, the proof of the following theorem (see [7, p. 224]) will not require knowledge of meromorphic continuation of Eisenstein series ($\mathcal{H} = \{z \mid \mathrm{Im}\, z > 0\}$):

THEOREM 16.4. *Let $\Gamma_0 \subseteq G = \mathrm{SL}(2, \mathbf{R})$ be a lattice such that $-I \in \Gamma_0$. Assume $\Lambda_0 = \left\{ \begin{pmatrix} 1 & n \\ 0 & 1 \end{pmatrix} \mid n \in \mathbf{Z} \right\}$ is a maximal unipotent subgroup of $\Gamma_0$. The Eisenstein series*

$$(16.5) \qquad E(z, s) = \frac{1}{2} \sum_{\gamma \in \Gamma_0 / \Lambda_0} (\mathrm{Im}\, \gamma^{-1} z)^s \qquad (z \in \mathcal{H},\ \mathrm{Re}\, s > 1)$$

*is convergent for* $\operatorname{Re} s > 1$. *Moreover, if* $U(\sigma)$, $0 < \sigma < \infty$, *is the set*

$$U(\sigma) = \{s \mid \operatorname{Re} s > 1, \frac{|s-1|}{\operatorname{Re} s - 1} < \sigma\},$$

*then for all* $\sigma$

(16.6) $$\lim_{\substack{s \to 1 \\ s \in U(\sigma)}} (s-1)E(z,s) = |\Gamma_0 \backslash \mathcal{H}|^{-1}$$

*where* $|\cdot|$ *denotes Poincaré volume.*

To relate Theorems 16.1 and 16.4 let $\Gamma$, $\Lambda$ and $v$ be as in the statement of Theorem 16.1. Choose $g_0 \in G$ so that $g_0^{-1} \Lambda g_0 = \Lambda_0 = \left\{ \begin{pmatrix} 1 & n \\ 0 & 1 \end{pmatrix} \mid n \in \mathbf{Z} \right\}$, and define $\Gamma_0 = g_0^{-1} \Gamma g_0$. Let $v_0 = \begin{pmatrix} 1 \\ 0 \end{pmatrix}$. Then replace $g_0$ by $-g_0$, if necessary, and reletter so that $g_0 v_0 = tv$ for some $t > 0$. Clearly,

(16.7) $$\operatorname{Card}(g\Gamma v \cap B(0,R)) = \operatorname{Card}(gg_0 \Gamma_0 v_0 \cap B(0,tR)).$$

Therefore, $c(\Gamma, v) = t^2 c(\Gamma_0, v_0)$, or since $|\Gamma_0 \backslash \mathcal{H}| = |\Gamma \backslash \mathcal{H}|$,

(16.8) $$c(\Gamma, v) = t^2 |\Gamma \backslash \mathcal{H}|^{-1}.$$

In order to adapt the present discussion to the requirements of Section 15 we require a lemma below. Note the identity

(16.9) $$(\operatorname{Im} g^{-1} i)^s = \|gv_0\|^{-2s} \qquad (g \in G, \ s \in \mathbf{C}).$$

LEMMA 16.10. *Let* $\Gamma$, $\Lambda$ *and* $v$ *be as in the statement of Theorem* 16.1. *There exists* $\tau = \tau(\Gamma, v) < \infty$ *such that*

(16.11) $$\operatorname{Card}(g\Gamma v \cap B(0,R)) < \tau(R^2 + 1) \qquad (R > 0, \ g \in G).$$

*Proof.* In view of (16.7) it is no loss of generality to suppose $\Gamma = \Gamma_0$, $\Lambda = \Lambda_0$ and $v = v_0$ (above). If $\mathcal{H}(R) = \{z \in \mathcal{H} \mid \operatorname{Im} z > \frac{1}{R^2}\}$, then (16.9) implies $g\gamma v_0 \in B(0,R)$ if, and only if, $\gamma^{-1} g^{-1} i \in \mathcal{H}(R)$. Since $\Lambda_0 v_0 = v_0$ and $\Lambda_0 \mathcal{H}(R) = \mathcal{H}(R)$, (16.11) is equivalent to a bound

(16.12) $$\operatorname{Card}\left(\Lambda_0 \backslash \left(\Gamma_0 g^{-1} i \cap \mathcal{H}(R)\right)\right) < \tau(R^2 + 1) \qquad (R > 0, \ g \in G).$$

Let $\mathcal{D}$ be a pairwise disjoint collection of open horodiscs such that: (a) $\Gamma_0 \mathcal{D} = \mathcal{D}$ and (b) if $|\mathcal{D}| = \bigcup_{D \in \mathcal{D}} D$, then $\Gamma_0 \backslash (\mathcal{H} \backslash |\mathcal{D}|)$ is compact. Observe that if $D \in \mathcal{D}$ and $\gamma \in \Gamma_0$ are such that $(\gamma D) \cap D \neq \emptyset$, then $\gamma D = D$. We shall divide (16.12) into two parts, one for $g^{-1} i \in |\mathcal{D}|$ and one for $g^{-1} i \in \mathcal{H} \backslash |\mathcal{D}|$.

Since $\Lambda_0 \backslash \mathcal{H}(R)$ has volume $R^2$ relative to the Poincaré volume $\frac{i}{2} \frac{dz \wedge d\bar{z}}{(\operatorname{Im} z)^2}$, and since $\Gamma_0 \backslash (\mathcal{H} \backslash |\mathcal{D}|)$ is compact, there exists $\tau_1 < \infty$ such that (16.12) is true with $\tau_1$ in place of $\tau$ and $g^{-1} i \notin |\mathcal{D}|$ in place of $g \in G$.

Let $\mathcal{H}_{0,1}(R) = \{z \in \mathcal{H}(R) \mid 0 \leq \operatorname{Re} z < 1\}$. Elementary euclidean geometry plus the fact the elements of $\mathcal{D}$ are pairwise disjoint imply there exists $\tau_2 < \infty$ such that $\operatorname{Card}\{D \in \mathcal{D} \mid D \cap \mathcal{H}_{0,1}(R) \neq \emptyset\} < \tau_2(R^2 + 1)$, $R > 0$. (Since $D_0 \cap \mathcal{H}_{0,1}(R) \neq \emptyset$ for all $R > 0$, it is necessary to use $R^2 + 1$ instead of $R^2$.) Now (16.12) is also true with $\tau_2$ in place of $\tau$ and $g^{-1}i \in |\mathcal{D}|$ in place of $g \in G$. Set $\tau = \operatorname{Max}(\tau_1, \tau_2)$, and (16.12) follows. The lemma is proved. □

*Proof of Theorem* 16.1. Let $\psi \in C_c(\mathbf{R}^2)$, and define $\hat{\psi}(g\Gamma) = \sum_{w \in \Gamma v} \psi(gw)$. Lemma 16.10 implies $\hat{\psi}$ is uniformly bounded on $G/\Gamma$. Since $\hat{\psi}$ is also continuous on $G/\Gamma$, the Eskin-McMullen theorem and Theorems 5.19 and 10.1 may be applied as in Section 15 to establish the existence of $c(\Gamma, v)$. The theorem is proved. □

In the notation of Theorem 16.4 define $N(g, R) = \operatorname{Card}(g\Gamma_0 v_0 \cap B(0, R))$. Setting aside the issue of convergence, the Eisenstein series (16.5) may be represented for any $z = g^{-1}i$ and $s$, $\operatorname{Re} s > 1$ by

$$
(16.13) \qquad E(z, s) = \frac{1}{2} \sum_{\gamma \in \Gamma_0/\Lambda_0} (\operatorname{Im} \gamma^{-1} g^{-1} i)^s
$$
$$
= \frac{1}{2} \sum_{\gamma \in \Gamma_0/\Lambda_0} \|g\gamma v_0\|^{-2s}
$$
$$
= \frac{1}{2} \int_0^\infty \frac{dN(g, R)}{R^{2s}}.
$$

Define $R_0(g) > 0$ so that $g\Gamma_0 v_0 \cap B(0, 2R_0(g)) = \emptyset$. Treat (16.13) as an improper Stieltjes integral over $(R_0(g), \infty)$ and integrate by parts to find

$$
(16.14) \qquad (s-1)E(z, s) = s(s-1) \int_{R_0(g)}^\infty \frac{N(g, R)}{R^2} R^{1-2s} dR.
$$

The calculation is justified for $\operatorname{Re} s > 1$ by Lemma 16.10. Convergence of (16.5) for $\operatorname{Re} s > 1$ is now established.

LEMMA 16.15. *Let $U(\sigma)$, $0 < \sigma < \infty$ be as in the statement of Theorem* 16.4. *We have for all $z \in \mathcal{H}$ and $0 < \sigma < \infty$*

$$
(16.16) \qquad \lim_{\substack{s \to 1 \\ s \in U(\sigma)}} (s-1)E(z, s) = \frac{c(\Gamma_0, v_0)\pi}{2}.
$$

*Proof.* Theorem 16.1 implies $N(g, R) = (c(\Gamma_0, v_0)\pi + \delta(g, R)) R^2$, where $\lim_{R \to \infty} \delta(g, R) = 0$. Substitute in (16.14) to obtain that (16.16) holds provided

$$
(16.17) \qquad \lim_{\substack{s \to 1 \\ s \in U(\sigma)}} s(s-1) \int_{R_0(g)}^\infty \delta(g, R) R^{1-2s} dR = 0.
$$

Since $|R^{1-2s}| = R^{1-2\operatorname{Re}s}$, and since $\left|\frac{s-1}{\operatorname{Re}s-1}\right| < \sigma$, $s \in U(\sigma)$, (16.17) follows from $\lim_{R\to\infty} \delta(g,R) = 0$.

Let $\mathcal{H}_{0,1} = \mathcal{H}_{0,1}(\infty) = \{z \in \mathcal{H} \mid 0 \leq \operatorname{Re}z < 1\}$. Since $\mathcal{H}_{0,1} \cong \Lambda_0\backslash\mathcal{H}$, $\mathcal{H}_{0,1}$ contains a geodesically convex fundamental domain $\Omega$ for $\Gamma_0$. We may suppose there exists $r_0$ such that $D_0 \cap \mathcal{H}_{0,1} \subseteq \Omega$, where $D_0$ is a horodisc $D_0 = \{z \in \mathcal{H} \mid \operatorname{Im}z > r_0^2\}$. There exists a function $0 \leq y(x) < r_0^2$ on $[0,1)$ such that

(16.18) $$\Omega = \{z = x + iy \mid 0 \leq x < 1,\ y > y(x)\}.$$

Let $\Gamma_0^* = \Gamma_0\backslash\Lambda_0$. We have for all $(z, s)$ that

(16.19) $$E(z,s) - (\operatorname{Im}z)^s = \frac{1}{2}\sum_{\gamma\in\Gamma_0^*/\Lambda_0}(\operatorname{Im}\gamma^{-1}z)^s \qquad (z \in \mathcal{H},\ \operatorname{Re}s > 1).$$

We shall be interested in (16.19) for $z \in \Omega$. In this case we claim the series on the right has no term such that $\gamma^{-1}z \in D_0$. Indeed, $\Omega$ is a fundamental domain which contains $D_0 \cap \mathcal{H}_{0,1}$ meaning $\gamma\Omega \cap D_0 = \emptyset$, $\gamma \in \Gamma_0^*/\Lambda_0$.

Let $N^*(g, R)$ be the counting function for $\Gamma_0^* v_0$. We have for $z = g^{-1}i$

(16.20)
$$E(z,s) - (\operatorname{Im}z)^s = \lim_{R_1\to\infty}\frac{1}{2}\int_{1/r_0}^{R_1}\frac{dN^*(g,R)}{R^{2s}}$$
$$= \lim_{R_1\to\infty}\left[\frac{N^*(g,R_1)}{2R_1^{2s}} + s\int_{1/r_0}^{R_1}\frac{N^*(g,R)}{R^{2s+1}}dR\right].$$

For any fixed $s$, $\operatorname{Re}s > 1$, Lemma 16.10 implies this convergence is uniform on $\Omega$. For fixed $R_1$ the first integral in (16.20) may be expressed as

(16.21) $$I_{R_1}(z,s) = \frac{1}{2}\sum_{\gamma\in\Gamma_0^*/\Lambda_0}(\operatorname{Im}(\gamma^{-1}z))^s \chi_{\mathcal{H}(R_1)}(\gamma^{-1}z) \qquad (z \in \Omega,\ \operatorname{Re}s > 1).$$

Using $G$-invariance of the volume element $\frac{i}{2}\frac{dz\wedge d\bar{z}}{(\operatorname{Im}z)^2}$ we have from (16.21)

$$\int_\Omega I_{R_1}(z,s)\frac{i}{2}\frac{dz\wedge d\bar{z}}{(\operatorname{Im}z)^2} = \frac{1}{2}\sum_{\gamma\in\Gamma_0^*/\Lambda_0}\int_{\gamma^{-1}\Omega\cap\mathcal{H}(R_1)}(\operatorname{Im}z)^s\frac{i}{2}\frac{dz\wedge d\bar{z}}{(\operatorname{Im}z)^2}.$$

Since $\lim_{R_1\to\infty}I_{R_1}(z,s) = E(z,s) - (\operatorname{Im}z)^s$ boundedly, we have for $\operatorname{Re}s > 1$

(16.22)
$$\int_\Omega (E(z,s) - (\operatorname{Im}z)^s)\frac{i}{2}\frac{dz\wedge d\bar{z}}{(\operatorname{Im}z)^2} = \int_{\mathcal{H}_{0,1}\backslash\Omega}(\operatorname{Im}z)^{s-2}\frac{i}{2}dz\wedge d\bar{z}$$
$$= \int_0^1 \frac{y(x)^{s-1}}{s-1}dx.$$

Finally, we observe that

LEMMA 16.23. *Let $0 < \sigma < \infty$, and let $U(\sigma, 1) = \{s \mid \operatorname{Re} s > 1, \frac{|s-1|}{\operatorname{Re} s - 1} < \sigma$ and $|s-1| < 1\}$. The product $(s-1)(E(z,s) - (\operatorname{Im} z)^s)$ is uniformly bounded on $\Omega \times U(\sigma, 1)$.*

*Proof.* From (16.20) (after $R_1 \to \infty$) and Lemma 16.10 we have

$$|(s-1)(E(z,s) - (\operatorname{Im} z)^s)| \leq |s(s-1)|\tau \int_{1/r_0}^{\infty} \frac{(R^2 + 1)}{R^{2\operatorname{Re} s + 1}} dR$$

$$= O\left(\frac{|s(s-1)|}{\operatorname{Re} s - 1}\right) = O(\sigma).$$

The lemma follows. $\square$

*Proof of Theorem* 16.4. Fix $0 < \sigma < \infty$. Lemma 16.15 implies that for each $z \in \mathcal{H}$

$$\lim_{\substack{s \to 1 \\ s \in U(\sigma, 1)}} (s-1)\left(E(z,s) - (\operatorname{Im} z)^s\right) = \frac{c(\Gamma_0, v_0)\pi}{2}.$$

Lemma 16.23, the bounded convergence theorem and (16.22) imply

$$\frac{c(\Gamma_0, v_0)\pi}{2}|\Omega| = \lim_{\substack{s \to 1 \\ s \in U(\sigma, 1)}} \int_0^1 y(x)^{s-1} dx$$

$$= 1.$$

Therefore, $\frac{c(\Gamma_0, v_0)\pi}{2} = |\Omega|^{-1}$ and Theorem 16.4 is proved. $\square$

RICE UNIVERSITY, HOUSTON, TX
*E-mail address*: veech@rice.edu